\documentclass[11pt]{amsart}
\usepackage{amsxtra}
\usepackage{amssymb}
\usepackage[mathscr]{eucal}
\theoremstyle{plain}

\usepackage{amsthm}
\newtheorem{manualtheoreminner}{Theorem}
\newenvironment{manualtheorem}[1]{%
  \renewcommand\themanualtheoreminner{#1}%
  \manualtheoreminner
}{\endmanualtheoreminner}

\newtheorem{theorem}{Theorem}[subsection]
\newtheorem*{theorem*}{Theorem}
\newtheorem{lemma}[theorem]{Lemma}
\newtheorem{definition-theorem}[theorem]{Definition-Theorem}
\newtheorem{proposition}[theorem]{Proposition}
\newtheorem{corollary}[theorem]{Corollary}
\newtheorem{definition}[theorem]{Definition}
\newtheorem{example}[theorem]{Example}
\newtheorem{remark}[theorem]{Remark}
\newtheorem{conjecture}[theorem]{Conjecture}
\newtheorem{notation}[theorem]{Notation}

\newtheorem*{maintheorem*}{Main Theorem}
\newcommand \bth[1] { \begin{theorem}\label{t#1} }
\newcommand \ble[1] { \begin{lemma}\label{l#1} }

\newcommand \bpr[1] { \begin{proposition}\label{p#1} }
\newcommand \bco[1] { \begin{corollary}\label{c#1} }
\newcommand \bde[1] { \begin{definition}\label{d#1}\rm }
\newcommand \bex[1] { \begin{example}\label{e#1}\rm }
\newcommand \bre[1] { \begin{remark}\label{r#1}\rm }
\newcommand \bcj[1] { \begin{conjecture}\label{j#1}\rm }
\newcommand \bqu[1]  { \medskip\noindent{\it{Question #1}} }
\newcommand \bnota[1] { \begin{notation}\label{n#1}\rm }
\renewcommand {\eth} { \end{theorem} }
\newcommand {\ele} { \end{lemma} }

\newcommand {\epr} { \end{proposition} }
\newcommand {\eco} { \end{corollary} }
\newcommand {\ede} { \end{definition} }
\newcommand {\eex} { \end{example} }
\newcommand {\ere} { \end{remark} }
\newcommand {\ecj} { \end{conjecture} }
\newcommand {\equ} {\smallskip}
\newcommand {\enota} { \end{notation} }
\newcommand \thref[1]{Theorem \ref{t#1}}
\newcommand \leref[1]{Lemma \ref{l#1}}
\newcommand \prref[1]{Proposition \ref{p#1}}
\newcommand \coref[1]{Corollary \ref{c#1}}

\newcommand \exref[1]{Example \ref{e#1}}
\newcommand \reref[1]{Remark \ref{r#1}}

\makeatletter
\newsavebox{\@brx}
\newcommand{\llangle}[1][]{\savebox{\@brx}{\(\m@th{#1\langle}\)}%
  \mathopen{\copy\@brx\kern-0.5\wd\@brx\usebox{\@brx}}}
\newcommand{\rrangle}[1][]{\savebox{\@brx}{\(\m@th{#1\rangle}\)}%
  \mathclose{\copy\@brx\kern-0.5\wd\@brx\usebox{\@brx}}}
\makeatother

\usepackage{tikz}
\usetikzlibrary{matrix}
\usepackage{graphicx,ctable,booktabs}
\usetikzlibrary{arrows.meta}
\usepackage{tikz-cd}

\usepackage{hyperref}

\DeclareMathOperator{\PHom}{PHom}

\DeclareMathOperator{\Ext}{Ext}

\DeclareMathOperator{\Spec}{Spec}

 \DeclareMathOperator{\Proj}{Proj}

 \DeclareMathOperator{\Hom}{Hom}

 \DeclareMathOperator{\im}{im}

\DeclareMathOperator{\Ind}{Ind}
\DeclareMathOperator{\ThickId}{ThickId}
\DeclareMathOperator{\Mod}{{\sf Mod}}
\DeclareMathOperator{\modd}{{\sf mod}}
\DeclareMathOperator{\stmod}{{\sf stmod}}
\DeclareMathOperator{\st}{{\sf st}}
\DeclareMathOperator{\St}{{\sf St}}
\DeclareMathOperator{\drin}{{\sf Z}}
\DeclareMathOperator{\StMod}{{ \sf StMod}}

\DeclareMathOperator{\Spc}{Spc}

\DeclareMathOperator{\cop}{cop}
\DeclareMathOperator{\op}{op}

\DeclareMathOperator{\cl}{cl}
\DeclareMathOperator{\Sp}{sp}

\newcommand{\mf}{\mathfrak}
\newcommand{\mc}{\mathcal}

\newcommand{\id}{\operatorname{id}}

\newcommand{\kk}{\Bbbk}
\newcommand{\bT}{\mathbf T}

\newcommand{\bS}{\mathbf S}
\newcommand{\bC}{\mathbf C}
\newcommand{\bK}{\mathbf K}

\newcommand{\bP}{\mathbf P}
\newcommand{\bQ}{\mathbf Q}
\newcommand{\bI}{\mathbf I}
\newcommand{\bJ}{\mathbf J}
\newcommand{\bD}{\mathbf D}

\newcommand{\Loc}{\operatorname{Loc}}
\newcommand{\XX}{\mathcal X}

\newcommand{\unit}{\ensuremath{\mathbf 1}}

\newcommand{\cohom}{\operatorname{H}^\bullet}

\newcounter{listequation}

\numberwithin{equation}{subsection}

\begin{document}

\title{Balmer Spectra and Drinfeld Centers}

\author[Kent B. Vashaw]{Kent B. Vashaw}
\address{
Department of Mathematics\\
Massachusetts Institute of Technology\\
Cambridge, MA 02139\\
U.S.A.}
\thanks{Research of K.B.V. was supported in part by a Board of Regents LSU fellowship, an Arthur K. Barton Superior Graduate Student Scholarship in Mathematics from LSU, NSF grant DMS-1901830, and NSF Postdoctoral Fellowship DMS-2103272.}
\email{kentv@mit.edu}

\begin{abstract}
The Balmer spectrum of a monoidal triangulated category is an important geometric construction which is closely related to the problem of classifying thick tensor ideals. We prove that the forgetful functor from the Drinfeld center of a finite tensor category ${\mathbf{C}}$ to ${\mathbf{C}}$ extends to a monoidal triangulated functor between their corresponding stable categories, and induces a continuous map between their Balmer spectra. 
We give conditions under which it is injective, surjective, or a homeomorphism. We apply this general theory to prove that Balmer spectra associated to finite-dimensional cosemisimple quasitriangular Hopf algebras (in particular, group algebras in characteristic dividing the order of the group) coincide with the Balmer spectra associated to their Drinfeld doubles, and that the thick ideals of both categories are in bijection. An analogous theorem is proven for certain Benson--Witherspoon smash coproduct Hopf algebras, which are not quasitriangular in general. 
\end{abstract}
\maketitle
\section*{Introduction}
\label{s:Intro}

Tensor triangular geometry, initiated by Balmer in \cite{Balmer1, Balmer2}, has proven to be a useful prism through which modular representation theory, algebraic geometry, commutative algebra, algebraic topology, and homotopy theory may all  be studied (for a few examples, see \cite{BKN1, BKN2, MT1, Balmer3, BS1}). The uniting feature is the existence, in each case, of a braided monoidal triangulated category; the braiding condition implies that there is a natural isomorphism
\[
X \otimes Y \cong Y \otimes X
\]
for all objects $X$ and $Y$. A noncommutative analogue of Balmer's theory (that is, one with no assumption of a braiding) was initiated and explored in \cite{BKS,NVY1, NVY2}, motivated by the abundance of examples of non-braided monoidal triangulated categories arising in representation theory. This theory defines a topological space, called the {\it Balmer spectrum}, for any monoidal triangulated category $\bT$. This space is denoted $\Spc \bT$, and is defined as the collection of prime ideals of $\bT$, reflecting the usual notion of prime spectrum from ring theory.

Non-braided monoidal triangulated categories arise naturally as the stable categories of finite tensor categories. Broadly speaking, if $\bC$ is a finite tensor category, then the stable category of $\bC$, denoted $\st(\bC)$, is the category obtained by factoring out the projective objects of $\bC$. One motivation for factoring out projectives comes from the theory of support varieties, where the support variety of an object only distinguishes an object up to direct sums with projective objects. The stable category is a monoidal triangulated category, where the monoidal product of $\st(\bC)$ is an extension of the monoidal product of $\bC$. 

An important tool in the study of tensor categories is the Drinfeld center, a categorical analogue of the center of a ring; it is a generalization of the quantum or Drinfeld double construction for Hopf algebras, originally introduced by Drinfeld in \cite{Drinfeld}. For any tensor category $\bC$, its Drinfeld center $\drin(\bC)$ is a braided tensor category equipped with a functor $F: \drin(\bC) \to \bC$. The Drinfeld center satisfies the universal property: if $G: \bD \to \bC$ is a strict tensor functor between strict tensor categories, and $\bD$ is braided, such that $G$ is bijective on objects and surjective on morphisms, then there exists a strict tensor functor $H: \bD \to \drin(\bC)$
\begin{center}
\begin{tikzcd}
\bC                       &                                             \\
\drin(\bC) \arrow[u, "F"] & \bD \arrow[lu, "G"'] \arrow[l, "H", dashed]
\end{tikzcd}
\end{center}
with $F \circ H = G$. 

If $\bC$ is abelian, then $\drin(\bC)$ is automatically abelian as well. We do not see an analogue for this argument in the triangulated case: if $\bT$ is triangulated, it does not seem to follow immediately that $\drin(\bT)$ is triangulated. This is a reflection of the fact that the morphism given in the extension axiom for triangulated categories is not necessarily unique. 

However, if $\bC$ is a finite tensor category, one can form its stable category $\st(\bC)$ on one hand; on the other hand, $\drin(\bC)$ is again a finite tensor category, and one can form its stable category $\st(\drin(\bC))$. The natural question that arises is, therefore: how are the Balmer spectra between these two categories connected?

This question is of particular interest because Balmer spectra are related intimately with cohomological support varieties (as in \cite{BPW1}); for example, the projectivization of the spectrum of the cohomology ring of the small quantum groups $u_{\zeta} (\mf{b})$ of Borel subalgebras at roots of unity (as computed in \cite{GK1, BNPP1}) identifies with the Balmer spectrum of its stable category \cite{NVY1}, which can be used to show that the support varieties for the small quantum Borel possess the tensor product property \cite{NP1, NVY2}. In many specific cases, for instance see \cite{FN1, Negron1, NP2}, the cohomology of Drinfeld doubles has been studied, and its relationship to the cohomology of the original finite tensor category explored. 

Additionally, this project will provide tools to aid in thick ideal classification problems. Balmer spectra, which are defined as the collection of prime ideals of the category, are intimately related to these problems, since every thick ideal of a rigid monoidal triangulated category is equal to an intersection of prime ideals. Classifications of thick ideals in various settings have been undertaken in many different settings, for instance in various categories arising from
\begin{enumerate}
\item commutative algebra and algebraic geometry \cite{Hopkins1,Thomason1, MT1};
\item Lie superalgebras \cite{BKN1};
\item finite groups and finite group schemes \cite{BCR1, FP1};
\item tilting modules for quantum groups and algebraic groups in positive characteristic \cite{Ostrik1, AHR1};
\item Hopf algebras which are not necessarily commutative, cocommutative, or even quasitriangular \cite{BW1, BKN1, NVY1, NVY2}.
\end{enumerate}
There are examples, for instance the small quantum groups of Borel subalgebras $u_{\zeta} (\mf{b})$ at roots of unity, where, as mentioned above, the Balmer spectrum and thick ideals are known for its stable module category; however, it is an open question to classify the Balmer spectrum and thick ideals for the stable category of its Drinfeld center, that is, the stable module category of $u_{\zeta} (\mf{g}) \otimes \kk T$, where $\kk T$ is the group algebra of the group of generators $K_i$ for $u_{\zeta} ( \mf{g})$. This motivates our central question, to reiterate: what relationship exists between the Balmer spectra of $\st(\bC)$ and $\st(\drin(\bC))$?

We answer this question by the following approach. 

In Section \ref{s:Prelim}, we give a brief background on tensor triangular geometry, compactly-generated triangulated categories, stable categories and finite tensor categories, support data, and Drinfeld centers, and establish notation.

Next, in Section \ref{s:DCBS}, we consider directly the relationship between the Balmer spectra $\Spc  \st(\bC)$ and $\Spc \st(\drin(\bC))$. Since the prime ideals of the Balmer spectrum of a non-braided monoidal triangulated category are a categorical analogue of the prime ideals in a noncommutative ring, we are motivated by prime ideal contraction, that is, the statement that if $\mf{p}$ is a prime ideal of a noncommutative ring $R$, then $\mf{p} \cap Z(R)$ is a prime ideal in $Z(R)$, the center of $R$. For general background on prime ideals for noncommutative rings, see \cite[Chapter 3]{GW1}. Finding a categorical analogue to prime ideal contraction is complicated by the fact that we work with $\st(\drin(\bC))$ rather than $\drin(\st(\bC))$; the latter is equipped with a forgetful functor $\drin(\st(\bC)) \to \st(\bC)$, but, as noted above, is not necessarily triangulated. Nevertheless, we verify that the forgetful functor $F: \drin(\bC) \to \bC$ extends to a functor $\overline{F}: \st(\drin(\bC)) \to \st(\bC)$. 

Reflecting the analogous property for rings, Balmer spectra of braided monoidal triangulated categories are functorial; but in the non-braided situation, a monoidal triangulated functor does not necessarily induce a continuous map between Balmer spectra. However, we show that $\overline{F}$ does induce a continuous map, and we obtain an analogue of prime ideal contraction. This is summarized by the following: 

\begin{manualtheorem}{A}
\label{intro-prime-id}{(See \prref{funct-drinf-cent-st} and \prref{fbar-induces-spc}).}
Let $\bC$ be a finite tensor category. There exists a monoidal triangulated functor $\overline{F}: \st(\drin(\bC)) \to \st(\bC)$ extending the forgetful functor $F: \drin(\bC) \to \bC$, which induces a continuous map $f: \Spc \st( \bC) \to \Spc \st(\drin(\bC))$, defined by 
$$f: \bP \mapsto \{ X \in  \st(\drin(\bC)) : \overline{F}(X) \in \bP\},$$
for $\bP \in \Spc \st(\bC)$.
\end{manualtheorem}

To study the image of the map $f$, we utilize the machinery of localization and colocalization functors. To apply these functors, one must work in the setting of compactly-generated triangulated categories. For us, the role of compactly-generated monoidal triangulated category will be filled by the stable category of the indization of $\bC$; this category will be referred to as $\St(\bC)$, and it contains $\st(\bC)$ as a triangulated subcategory. For the details of this setting, see Section \ref{ss:compact-gen}. It is straightforward that the functor $\overline{F}$ extends to a functor $\St(\drin(\bC)) \to \St(\bC)$; denote this extension again by $\overline{F}$. We are then able to use the kernel of this functor to describe the image of $f$.

\begin{manualtheorem}{B}
\label{intro-im-f}{(See \prref{map-f-surj}).}
Denote by $\bK$ the kernel of $\overline{F}: \St(\drin(\bC)) \to \St(\bC)$, and  $f: \Spc \st( \bC) \to \Spc \st(\drin(\bC))$ the continuous map induced by $\overline{F}$ as above. Then there are containments
\begin{align*}
\{  \bP \in \Spc \st(\bC) : \bP \supseteq \bK \cap \st(\drin(\bC)) \} &\supseteq \im f\\
& \supseteq \{  \bP \in \Spc \st(\drin(\bC)) : \Loc(\bP) \supseteq \bK \}.\end{align*}
Here, $\Loc(\bP)$ refers to the localizing subcategory (meaning triangulated and closed under set-indexed coproducts) of $\St(\drin(\bC))$ generated by $\bP$.
\end{manualtheorem}

This implies that if $\bC$ satisfies the following property, then $f$ is surjective.
\begin{align*} 
\mbox{Property (*): For $X$ in $\Ind(\drin(\bC))$, if $F(X)$ is projective, then so is $X$. \label{EE:propp} } 
\end{align*} 
Additionally, if $\bC$ is a braided tensor category to begin with, then we prove that $f$ is injective. This leads to the following theorem.

\begin{manualtheorem}{C}
\label{intro-surj-inj}{(See \thref{f-inj-surj-homeo}).}
Let $\bC$ be a finite braided tensor category satisfying Property (*). Then $f$ is a homeomorphism $\Spc \st(\bC) \xrightarrow{\cong} \Spc \st(\drin(\bC))$, and there is a bijection between the thick ideals of $\st(\drin(\bC))$ and the thick ideals of $\st(\bC)$, given by  
$$\bI \mapsto \langle \overline{F}(X) : X \in \bI \rangle$$ for a thick ideal $\bI$ of $\st(\drin(\bC))$.
\end{manualtheorem}

In Section \ref{s:applications}, we illustrate the theory with concrete examples. We first consider $\bC$ to be alternately $\modd(\kk G)$ and $\modd((\kk[G]^{\cop})$, for $G$ a finite group and $\kk$ an algebraically closed field of characteristic $p$ dividing the order of $G$, where $\kk G$ denotes the group algebra of $G$ and $\kk[G]$ denotes the dual group algebra to $\kk G$. Of these two examples, the first satisfies Property (*) and the second does not. This allows us to classify the Balmer spectrum and classify the thick ideals for $\stmod ( D(\kk G))$, where $D(\kk G)$ is the Drinfeld double of the group algebra $\kk G$. We are then able to generalize this example in the following way. 

\begin{manualtheorem}{D}
\label{intro-class-thm}
{(See \prref{group-dualgroup}, \prref{quasi-cosemi}, and \thref{cosemi-bw}).}
For the following classes of Hopf algebras $H$, the Balmer spectrum of $\stmod(D(H))$ is homeomorphic via the map $f$ to the Balmer spectrum of $\stmod(H)$, and the thick ideals of the two categories are in bijection:
\begin{enumerate}
\item finite-dimensional cosemisimple quasitriangular Hopf algebras (e.g. group algebras of finite groups $G$ in characteristic dividing the order of $G$);
\item Benson--Witherspoon smash coproducts $(\kk[G] \# \kk L)^*$, where $G$ and $L$ are finite groups with dual group algebra and group algebra $\kk[G]$ and $\kk L$ respectively, $\kk$ an algebraically closed field of characteristic $p$ dividing the order of $G$ and not dividing the order of $L$, such that $L$ acts by group automorphisms on $G$.
\end{enumerate}
\end{manualtheorem}

\subsection*{Acknowledgements.} The author would like to thank Daniel Nakano and Milen Yakimov for many useful discussions, and also to thank Pavel Etingof, Siu-Hung Ng, Victor Ostrik, and Sean Sanford for providing valuable comments used to improve this paper. We also thank the anonymous referee for their careful reading of the paper and helpful suggestions.

\section{Preliminaries}
\label{s:Prelim}

\subsection{Tensor triangular geometry}
\label{ss:ttg}

We will recall some of the background of noncommutative tensor triangular geometry. Following the terminology of \cite{NVY1, NVY2}, a {\it monoidal triangulated category} $\bT$ is a category such that the following conditions hold:
\begin{enumerate}
\item $\bT$ is triangulated: it is an additive category equipped with an additive autoequivalence $\Sigma: \bT \to \bT$, called the {\it shift functor}, and a collection of distinguished triangles
$$A \to B \to C \to \Sigma A,$$ subject to the usual axioms (see \cite{Happel1,Neeman1});
\item $\bT$ is monoidal: it is equipped with a monoidal product $\otimes$ and unit $\unit_{\bT}$, subject to the usual associativity and unit axioms (see \cite{Kassel1, BK1, EGNO1});
\item the triangulated and monoidal structures on $\bT$ are compatible: for any object $A$ of $\bT$, the functors $A \otimes - $ and $- \otimes A$ are triangulated functors. In other words, there exists a natural isomorphism $\Sigma(A) \otimes B \cong \Sigma ( A \otimes B) \cong A \otimes \Sigma(B)$, such that if 
$$A \to B \to C \to \Sigma A$$
is a distinguished triangle, then for any object $D$, the triangles
$$D \otimes A \to D \otimes B \to D \otimes C \to \Sigma ( D \otimes A)$$
and
$$A \otimes D \to B \otimes D \to C \otimes D \to \Sigma( A\otimes D)$$
are distinguished. 
\end{enumerate}

\bre{tensor-tri}
In the terminology of \cite{Balmer1, Balmer2}, a {\it tensor triangulated category} is a monoidal triangulated category such that the monoidal product is symmetric. Note that contrary to the definition of tensor category as in \cite{EGNO1}, a tensor triangulated category is not required to have duals.
\ere

We will recall the definition of the Balmer spectrum of a monoidal triangulated category $\bT$, as in \cite{BKS,NVY1}.
\begin{enumerate}
\item A {\it (two-sided) thick ideal} $\bI$ of $\bT$ is a full subcategory such that the following hold.
	\begin{enumerate}
	\item $\bI$ is triangulated: it is closed under $\Sigma$ and $\Sigma^{-1}$, and if $$A \to B \to C \to \Sigma A$$ is a distinguished triangle, then if any two of $A, B,$ and $C$ are in $\bI$, then so is the third;
	\item $\bI$ is thick: if $A \oplus B$ is in $\bI$, then so are $A$ and $B$;
	\item $\bI$ is an ideal: if $A \in \bI$, then so are $A \otimes B$ and $B \otimes A$, for any object $B$. 
	\end{enumerate}
The collection of thick ideals of $\bT$ will be denoted by $\ThickId(\bT)$, and the thick ideal generated by a collection of objects $\mc{T}$ will be denoted $\langle \mc{T} \rangle$.
\item A thick ideal $\bP$ of $\bT$ is called {\it prime} if for all thick ideals $\bI$ and $\bJ$ of $\bT$, a containment $\bI \otimes \bJ \subseteq \bP$ implies either $\bI$ or $\bJ \subseteq \bP$; equivalently, $\bP$ is prime if and only if for all objects $A$ and $B$ of $\bT$, a containment $A \otimes \bT \otimes B \subseteq \bP$ implies either $A$ or $B$ is in $\bP$ (see \cite[Theorem 3.2.2]{NVY1}). Here $\bI \otimes \bJ$ refers to the collection of objects $\{A \otimes B  : A \in \bI, \; B \in \bJ\}$, and $A \otimes \bT \otimes B$ refers to the collection of objects $\{ A \otimes C \otimes B  : C \in \bT \}$ for $A$ and $B$ in $\bT$.
\item A thick ideal $\bP$ of $\bT$ is called {\it completely prime} if $A\otimes B \in \bP$ implies either $A$ or $B \in \bP$, for all objects $A$ and $B$ of $\bT$.
\item The {\it Balmer spectrum} of $\bT$, denoted $\Spc \bT$, is the collection of prime ideals of $\bT$ under the Zariski topology, where closed sets are defined as the sets
$$V_{\bT}(\mc{T}) = \{ \bP \in \Spc \bT : \mc{T} \cap \bP = \varnothing \}$$
for all collections of objects $\mc{T}$ of $\bT$.  
\item An arbitrary open set of $\Spc \bT$, that is, the complement of a closed set $V_{\bT}(\mc{T})$ for some collection $\mc{T}$ of objects of $\bT$, will be denoted 
\[
U_{\bT}(\mc{T}):=\Spc \bT \backslash V_{\bT}(\mc{T}) = \{ \bP \in \Spc \bT : \mc{T} \cap \bP \not = \varnothing\}.
\]
\end{enumerate}

Note that every completely prime ideal is prime. If $\bT$ is a braided category then every prime ideal is completely prime, and so in that case the two notions coincide. 

\bre{topol}
We emphasize that this choice of topology on the Balmer spectrum does not match what one might expect, by the analogy to ring theory. This reflects the fact that in natural examples when the Balmer spectrum of a monoidal triangular category $\bT$ is realized concretely as the $\Proj$ or $\Spec$ of a commutative ring $R$, the bijection between prime ideals of $\bT$ and the (homogeneous) prime ideals of $R$ is containment-reversing. See \exref{fin-group} below for concrete examples.
\ere

\bre{sheaf}
While we have only defined the Balmer spectrum as a topological space, Balmer's original definition \cite[Section 6]{Balmer1} gives $\Spc$ the additional structure of a ringed space (which is actually locally ringed, by \cite[Corollary 6.6]{Balmer2}). Many of the classification theorems for Balmer spectra prove existence of isomorphisms of ringed spaces, rather than just homeomorphisms of topological spaces. However, the ringed space structures will not play a role in this paper, so we omit the precise definition.
\ere

We recall one topological property of the Balmer spectrum, for reference. This was proven by Balmer \cite[Corollary 2.17]{Balmer1}.
\bth{balmer-top}
Let $\bT$ be a braided monoidal triangulated category. Then $\Spc \bT$ is Noetherian if and only if every closed subset of $\Spc \bT$ is of the form $V_{\bT}(A)$, for some object $A$ of $\bT$.
\eth

\bre{hochster}
In fact, if $\bT$ is braided, or if there is an object of $\bT$ which generates $\bT$ as a thick subcategory, then $\Spc \bT$ is a spectral topological space. In other words, $\Spc \bT$ is $T_0$, quasicompact, the quasicompact open sets form an open basis, and every non-empty irreducible closed subset has a generic point. It is a theorem of Hochster that this implies $\Spc \bT$ is homeomorphic to the prime spectrum of a commutative ring \cite[Theorem 6]{Hochster}. 
\ere

Suppose $\bT$ is rigid, in other words, every object is dualizable (cf. \cite[Section 2.10]{EGNO1}). We then obtain the following facts, which we recall for reference. Both follow directly from the fact that if $A$ is dualizable with dual $A^*$, then $A$ is a direct summand of $A \otimes A^* \otimes A$. 

\bpr{semiprime}
Let $\bT$ be a rigid monoidal triangulated category. Let $A$ be an object of $\bT$ with dual $A^*$. Then
\begin{enumerate}
\item $\langle A \rangle = \langle A^* \rangle$, and
\item every thick two-sided ideal $\bI$ of $\bT$ is semiprime, i.e. it is the intersection $$\bI = \bigcap_{\bP \in \Spc \bT, \; \bI \subseteq \bP} \bP$$ of prime ideals over itself. Equivalently, for every ideal $\bI$ of $\bT$, if the set of objects $A \otimes \bT \otimes A \subseteq \bI$ for some object $A$ in $\bT$, then $A \in \bI$, where $A \otimes \bT \otimes A$ refers to the collection $\{A \otimes B \otimes A : B \in \bT\}$. 
\end{enumerate}
\epr

For the details of the proofs, see \cite[Lemma 5.1.1]{NVY1} and \cite[Proposition 4.1.1]{NVY2}. 

\subsection{Compactly-generated triangulated categories}
\label{ss:compact-gen}

A powerful result in the theory of triangulated categories is Brown representability, which ensures the existence of adjoints to certain triangulated functors \cite[Chapter 8]{Neeman1}. However, in order to apply these results, one must work in the setting of compactly-generated triangulated categories. We recall the definition now.
\begin{enumerate}
\item An object $C$ in a triangulated category $\bT$ is {\it compact} if the functor $\Hom_{\bT}(C,-)$ commutes with arbitrary set-indexed coproducts. If $\bT$ is a triangulated category, then $\bT^c$ will denote the subcategory of compact objects.
\item A {\it localizing subcategory} of a triangulated category is a triangulated subcategory which is also closed under taking set-indexed coproducts. The smallest localizing category containing a collection $\mc{T}$ of objects will be denoted $\Loc(\mc{T})$ and will be referred to as the {\em{localizing category generated by $\mc{T}$}}.
\item A {\it compactly-generated triangulated category} is a triangulated category $\bT$ which contains arbitrary set-indexed coproducts such that $\Loc(\bT^c)= \bT$.
\end{enumerate}
Note that any localizing subcategory $\bI$ of $\bT$ is thick by a version of the Eilenberg swindle: if $A \oplus B$ is in $\bI$, then we have a distinguished triangle
\[
(A \oplus B)^{\oplus \mathbb{N}} \to (A \oplus B)^{\oplus \mathbb{N}} \to A \to \Sigma (A \oplus B)^{\oplus \mathbb{N}}, 
\]
where the first map sends the $i$th copy of $B$ in the first object to the $i$th copy of $B$ in the second object, and sends the $i$th copy of $A$ in the first object to the $(i+1)$th copy of $A$ in the second object. Since $\bI$ is localizing, the first and second objects are in $\bI$, and hence $A$ is in $\bI$ as well. For additional background on compactly-generated triangulated categories, see \cite[Section 1.3.9]{BIK2}.

The following theorem, due to Rickard \cite{Rickard1}, is the primary technical reason we need to move to the compactly-generated setting. For details, see \cite[Theorem 3.1.1, Theorem 3.1.2]{BKN1}, \cite[Section 3]{BIK1}, and \cite[Section 2]{BIK2}. 

\bth{local}
Let $\bT$ be a compactly generated triangulated category. Given a thick subcategory $\bS$ of $\bT^c$, there exist functors ${\it{\Gamma}}_{\bS}$ and $L_{\bS}$ from $\bT \to \bT$, which gives for every object $M$ of $\bT$ a distinguished triangle
$${\it{\Gamma}}_{\bS} (M) \to M \to L_{\bS}(M) \to \Sigma ({\it{\Gamma}}_{\bS}(M)),$$
such that
\begin{enumerate}
\item $L_{\bS}$ and ${\it{\Gamma}}_{\bS}$ are unique up to isomorphism,
\item ${\it{\Gamma}}_{\bS}(M)$ is in $\Loc(\bS)$,
\item $L_{\bS}(M)$ is in $\Loc(\bS)^\perp$, that is, there are no nonzero maps from $\Loc(\bS) \to L_{\bS} (M)$, and
\item $M \in \Loc(\bS)$ if and only if ${\it{\Gamma}}_{\bS} (M) \cong M$, or, equivalently, $L_{\bS} (M) \cong 0$.   
\end{enumerate}
\eth

The functors ${\it{\Gamma}}_{\bS}$ and $L_{\bS}$ are called {\it colocalizing} and {\it localizing} functors, respectively. They are constructed by first taking a Verdier quotient of $\bT$ by $\Loc(\bS)$, that is, forming a category where all morphisms with cones in $\Loc(\bS)$ are formally inverted, which one may do using the calculus of roofs. This quotient is a triangulated category where the objects isomorphic to 0 are precisely those from $\Loc(\bS)$, and Brown representability guarantees that there are right adjoint functors $i^!$ and $j_*$ to the inclusion $i_*: \Loc(\bS) \to \bT$ and quotient $j^*: \bT \to \bT/ \Loc(\bS)$ functors, giving a diagram
\[
\Loc(\bS) \begin{array}{c} {i_*} \atop {\longrightarrow} \\ {\longleftarrow}\atop{i^!} \end{array}   \bT \begin{array}{c} {j^*} \atop {\longrightarrow} \\ {\longleftarrow}\atop{j_*} \end{array}   \bT/ \Loc(\bS).
\]
The functor $L_{\bS}$ is then defined as $j_* \circ j^*$, and ${\it{\Gamma}}_{\bS}$ is defined as $i_* \circ i^!$. For the details of the categorical localization and Verdier quotient, as well as additional details on the formation of the localization and colocalization functors, see \cite[Section 2.1, Theorem 8.4.4]{Neeman1}, \cite{Krause1}, and \cite[Section 3]{Stevenson}.

\subsection{Stable categories and finite tensor categories}
\label{ss:stab}

The monoidal triangulated categories that are the primary focus of this paper arise as stable categories. We first recall the construction of the stable category of any quasi-Frobenius category. Recall that a quasi-Frobenius category is an abelian category with enough projectives, such that projective and injective objects coincide. For any quasi-Frobenius category $\bC$, one may form the stable category $\st (\bC)$, which is triangulated (see \cite[Chapter I]{Happel1}). The stable category is constructed by factoring out the projective objects of $\bC$. In more detail, let $\PHom_{\bC} (A,B)$ consist of the morphisms $f: A \to B$ in $\bC$ such that $f$ factors through a projective object. The stable category $\st(\bC)$ is the category where
\begin{enumerate}
\item objects are the same as the objects of $\bC$;
\item morphisms $A \to B$ are defined to be $\Hom_{\bC} (A, B) / \PHom_{\bC} (A, B).$
\end{enumerate}
There is a straightforward functor $G: \bC \to \st(\bC)$ sending objects to themselves and morphisms to their image in the quotient. 

If $P$ is a projective object of $\bC$, note that the corresponding object $G(P)$ in $\st (\bC)$ is isomorphic to 0, since $\id_P$ factors through a projective; and the converse is also true, since $G(P) \cong 0$ in $\st(\bC)$ implies that the 0-morphism $G(P) \to G(P)$ is equal to $\id_{G(P)}$ in $\Hom_{\st(\bC)} (G(P), G(P))$, in other words, $\id_P$ factors through a projective $Q$ in $\bC$: 
\begin{center}
\begin{tikzcd}
P \arrow[rd, hook] \arrow[rr, "\id_P"] &                         & P \\
                                       & Q \arrow[ru, two heads] &  
\end{tikzcd}
\end{center}
Of course, this implies $P$ is a summand of $Q$, and so $P$ is projective. 

We recall the triangulated structure on $\st(\bC)$, for reference. If $A$ is an object of $\bC$, denote by $\Omega(A)$ the kernel of the projective cover of $A$. The functor $\Omega$ extends to the stable module category, and this in fact gives an autoequivalence on $\st(\bC)$. The shift $\Sigma$ is then defined to be $\Sigma(A)= \Omega^{-1}(A)$. For any short exact sequence of objects in $\bC$, say
$$0 \to A \to B \to C \to 0,$$
there exists a triangle $$A \to B \to C \to \Sigma A$$ in $\st(\bC)$; the distinguished triangles of $\st(\bC)$ are then defined to be all triangles which are isomorphic to triangles of this form. 


We now specialize to the case that $\bC$ is a finite tensor category. Recall that a finite tensor category (following the notation given in \cite{EO1,EGNO1}) consists of a monoidal category $\bC$ such that:
\begin{enumerate}
\item $\bC$ is abelian and $\kk$-linear for an algebraically closed field $\kk$;
\item the tensor product $- \otimes -$ is bilinear on spaces of morphisms;
\item every object of $\bC$ has finite length;
\item $\Hom_{\bC} (\unit, \unit) \cong \kk$;
\item for any pair of objects $A$ and $B$, the vector space $\Hom_{\bC} (A, B)$ is finite-dimensional over $\kk$;
\item $\bC$ has enough projectives;
\item there are finitely many isomorphism classes of simple objects of $\bC$;
\item $\bC$ is rigid, i.e. every object has a left and a right dual. 
\end{enumerate}
The prototypical example of a finite tensor category is the category of finite-dimensional modules of a finite-dimensional Hopf algebra $H$. 
\begin{notation}
Denote the category of finite-dimensional modules of an algebra $H$ by $\modd(H)$. Denote the category of all (not necessarily finite-dimensional) modules of $H$ by $\Mod(H)$.
\end{notation}

Recall that if $\bC$ is a finite tensor category, it is a consequence that the tensor product is biexact (\cite[Proposition 4.2.1]{EGNO1}). Additionally, every finite tensor category is quasi-Frobenius \cite[Proposition 6.1.3]{EGNO1}. The stable category $\st(\bC)$ inherits a monoidal product directly from $\bC$: we define $G(A) \otimes G(B) := G ( A \otimes B)$, and similarly for morphisms $f: A \to B$ and $g: C \to D$ we define $G(f) \otimes G(g) := G(f \otimes g)$. This is well-defined: if $G(f) = G(\hat f)$, then $f-\hat f$ factors through a projective $P$, and then $f \otimes g - \hat f \otimes g$ factors through $P \otimes D$, which is projective by projectivity of $P$ (see \cite[Proposition 4.2.12]{EGNO1}).

Although the primary objects of focus in this paper are stable categories $\st(\bC)$ for finite tensor categories $\bC$, note that $\st(\bC)$ is {\it not} compactly-generated, since in particular it does not contain arbitrary set-indexed coproducts. Thus, in order to apply \thref{local}, it is necessary to produce a compactly-generated monoidal triangulated category which contains $\st(\bC)$ as a monoidal triangulated subcategory. In fact, this is possible, using the $\Ind$-completion (see \cite[Chapter 6]{KS1}) of $\bC$:
\bth{big-st}
Let $\bC$ be a finite tensor category. Then its $\Ind$-completion $\Ind(\bC)$ is a quasi-Frobenius abelian monoidal category, and its stable category $\st(\Ind(\bC))$ is a compactly-generated monoidal triangulated category, with $\st(\Ind(\bC))^c \cong \st(\bC)$ via the stabilization of the natural inclusion functor $\bC \to \Ind(\bC)$.
\eth

\begin{proof}
See \cite[Theorem A.0.1]{NVY3}.
\end{proof}

Concretely, the there exists a finite-dimensional algebra $A$ such that $\bC \cong \modd(A)$, the category of finite-dimensional $A$-modules \cite[page 10]{EGNO1}. Then $\Ind(\bC) \cong \Mod(A)$, the category of all $A$-modules. 

\begin{notation}
\label{big-st-name}
If $\bC$ is a finite tensor category, we denote 
\[
  \St(\bC):=\st(\Ind(\bC)),
\]
to avoid crowding the notation.
\end{notation}

\subsection{Support data}
\label{ss:supp}

Suppose that $\bT$ is a monoidal triangulated category and $S$ a topological space. We will denote the collection of subsets of $S$ by $\mc{X}(S)$, closed subsets of $S$ by $\mc{X}_{\cl}(S)$, and specialization-closed subsets of $S$ by $\mc{X}_{\Sp}(S)$; recall that by definition, a set is specialization-closed if it is a union of closed sets. When the underlying space is clear from context, we will denote these collections by $\mc{X}$, $\mc{X}_{\cl}$, and $\mc{X}_{\Sp}$. 

Given a monoidal triangulated category $\bT$ and a topological space $S$, a {\it support datum on $\bT$ with value in $S$} is a map $\sigma: \bT \to \mc{X}_{\cl}(S)$ satisfying the following axioms:
\begin{enumerate}
\item $\sigma(0)=\varnothing $ and $ \sigma(\unit)= S$;
\item $\sigma(A \oplus B)=\sigma(A) \cup \sigma(B)$, $\forall A,B \in \bT$; 
\item $\sigma(\Sigma A)=\sigma(A)$, $\forall A \in \bT$; 
\item if $A \to B \to C \to \Sigma A$ is a distinguished triangle, then $\sigma(A) \subseteq \sigma(B) \cup \sigma(C)$;
\item $\bigcup_{C \in \bT} \sigma (A   \otimes C \otimes B) = \sigma (A) \cap \sigma(B)$, $\forall A, B \in \bT$. 
\end{enumerate}

See \cite[Section 4]{NVY1} for a more in-depth discussion of support data (although note that in that paper, a support datum is permitted to take value in $\mc{X}(S)$ rather than $\mc{X}_{\cl}(S)$). For any monoidal triangulated category $\bT$, the map $V_{\bT} (A)= \{ \bP \in \Spc \bT : A \not \in \bP \}$ defined above is a support datum $\bT \to  \mc{X}_{\cl}(\Spc \bT)$, since by definition, $V_{\bT}(A)$ is a closed set in $\Spc \bT$. We will refer to this support datum as the {\it Balmer support}. Indeed, the Balmer support satisfies a universal property in the category of support data, see \cite[Theorem 4.2.2]{NVY1}. 

\bth{univ-balm}
If $\sigma: \bT \to \mc{X}_{\cl}(S)$ is a support datum with value in $S$, then there exists a unique continuous map
\[
S \xrightarrow{f} \Spc \bT
\]
such that $\sigma(A)= f^{-1}(V(A))$ for all $A \in \bT$.
\eth

For any support datum $\sigma$, we have a map $$\Phi_\sigma (\mc{T}) := \bigcup_{A \in \mc{T}} \sigma (A),$$ where $\mc{T}$ is any collection of objects of $\bT$. If $\sigma$ takes values in $\mc{X}_{\cl}$, then by definition $\Phi_\sigma$ takes values in $\mc{X}_{\Sp}$. The map $\Phi_\sigma$ in fact only depends on thick ideals rather than arbitrary subsets, since by \cite[Lemma 4.3.2]{NVY1} we have $\Phi_\sigma(\mc{T}) = \Phi_\sigma ( \langle \mc{T} \rangle)$.

For a support datum $\sigma$, we have a second map $\Theta_\sigma : \mc{X}_{\Sp} \to \ThickId(\bT)$ defined by $$\Theta_\sigma(S') := \{ A \in \bT : \sigma(A) \subseteq S'\}$$ for any specialization-closed subset $S'$ of $S$. 
For any specialization closed set $S'$, the collection $\Theta_{\sigma}(S')$ is a thick ideal of $\bT$. Hence, we have the following collection of maps, given a support datum $\sigma$ on $\bT$:
$$
\ThickId(\bT) \begin{array}{c} {\Phi_{\sigma}} \atop {\longrightarrow} \\ {\longleftarrow}\atop{\Theta_\sigma} \end{array}  \XX_{\Sp}.
$$
Classifications of thick ideals are obtained in many cases (see \cite{Balmer1, Balmer2, BKN1, BKN2, NVY1, NVY2} for examples) by constructing a support datum for which these maps are bijective and inverse to each other. In that case, the support datum $\sigma$ is called {\it classifying}. For rigid braided monoidal triangulated categories $\bT$, the Balmer support $V_{\bT}$ is always classifying \cite[Theorem 4.10]{Balmer1}.

\bex{fin-group}
For a finite group scheme $G$, the {\it cohomological support} is the map $\stmod(G) \to \mc{X}_{\cl}\left (\Proj \cohom(G,\kk)\right)$ defined by
\[
M \mapsto \{ \mf{p} \in \Proj \cohom(G,\kk) : \mf{p} \text{ contains }I(M)\},
\]
where $I(M)$ is the annihilator of $\bigoplus_{i \geq 0}\Ext^i_{G}(M,M)$ in $\cohom(G,\kk):=\bigoplus_{i \geq 0} \Ext^i_{G}(\unit,\unit)$ under the action induced by the functor $M \otimes -$ \cite[Section 5.7]{Benson1}. Cohomological support is a support datum; the most nontrivial property is (5), referred to as the {\it tensor product property}, and was proven by Friedlander--Pevtsova \cite{FP1}. It is a theorem that for finite group schemes, the cohomological support is classifying, and the map $f: \Proj \cohom(G,\kk) \to \Spc \stmod(G)$ is a homeomorphism \cite{BCR1, FP1, Balmer1}. Cohomological support exists for arbitrary finite tensor categories \cite{BPW1}, but is not known to be classifying in general, see \cite[Conjecture E]{NVY3}.
\eex

\subsection{The Drinfeld center}
\label{ss:drinf}

Let $\bC$ be a strict monoidal category. Then the {\it Drinfeld center} or {\it center} of $\bC$, which we will denote by $\drin(\bC)$, is defined as the following braided monoidal category.
\begin{enumerate}
\item Objects are pairs $(A, \gamma)$ where $A$ is an object of $\bC$ and $\gamma$ is a natural isomorphism $\gamma_B: B \otimes A \xrightarrow{ \cong} A \otimes B$ for all $B \in \bC$, satisfying the diagram
\begin{center}
\begin{tikzcd}
B \otimes C \otimes A \arrow[r, "\id_B \otimes \gamma_C"] \arrow[rr, "\gamma_{B \otimes C}", bend left] & B \otimes A \otimes C \arrow[r, "\gamma_B \otimes \id_C"] & A \otimes B \otimes C
\end{tikzcd}
\end{center}
for all $B$ and $C$. Such a natural isomorphism $\gamma$ is called a {\it half-braiding} of $A$.
\item Morphisms $(A, \gamma) \to (A', \gamma')$ are morphisms $f: A \to A'$ such that for all $B$, the diagram
\begin{center}
\begin{tikzcd}
B \otimes A \arrow[r, "\id_B \otimes f"] \arrow[d, "\gamma_B"] & B \otimes A' \arrow[d, "\gamma'_B"] \\
A \otimes B \arrow[r, "f \otimes \id_B"]                       & A' \otimes B                       
\end{tikzcd}
\end{center} 
commutes. 
\item The monoidal product $(A, \gamma) \otimes (A', \gamma')$ is defined as $(A \otimes A', \tilde \gamma)$ where $\tilde \gamma$ is defined as 
\begin{center}
\begin{tikzcd}
B \otimes A \otimes A' \arrow[d, "\tilde \gamma_B"] \arrow[r, "\gamma_B \otimes \id_{A'}"] & A \otimes B \otimes A' \arrow[ld, "\id_A \otimes \gamma'_B"] \\
A \otimes A'\otimes B                                                                      &                                                             
\end{tikzcd}
\end{center}
\item The braiding $c_{(A, \gamma), (A', \gamma')}: (A, \gamma) \otimes (A', \gamma') \xrightarrow{\cong} (A', \gamma') \otimes (A, \gamma)$ is defined as $\gamma'_A$. The map $\gamma'_A$ being a valid map in $\drin(\bC)$ amounts to checking the commutativity of the diagram
\begin{center}
\begin{tikzcd}
B \otimes A \otimes A' \arrow[r, "\id_B \otimes \gamma_A'"] \arrow[d, "\gamma_B \otimes \id_{A'}"] & B \otimes A' \otimes A \arrow[d, "\gamma'_B \otimes \id_A"]   \\
A \otimes B \otimes A' \arrow[d, "\id_A \otimes \gamma_B'"] & A' \otimes B \otimes A \arrow[d, "\id_{A'} \otimes \gamma_B"] \\
A \otimes A' \otimes B \arrow[r, "\gamma_A' \otimes \id_B"]                                        & A' \otimes A \otimes B                                       
\end{tikzcd}
\end{center}
This diagram commutes by the naturality of $\gamma'$, since it can be rewritten, using the defining diagram for $\gamma'$, as
\begin{center}
\begin{tikzcd}
B \otimes A \otimes A' \arrow[r, "\gamma'_{B \otimes A}"] \arrow[d, "\gamma_B \otimes \id_{A'}"] & A' \otimes B \otimes A \arrow[d, "\id_{A'} \otimes \gamma_B"] \\
A \otimes B \otimes A' \arrow[r, "\gamma'_{A \otimes B}"]                                        & A' \otimes A \otimes B                                       
\end{tikzcd}
\end{center} 
\end{enumerate}

We will denote by $F: \drin(\bC) \to \bC$ the forgetful functor sending $(A, \gamma) \mapsto A$. 

If $H$ is a Hopf algebra and $\bC$ is the category of $H$-modules, it is well-known that the Drinfeld center $\drin(\bC)$ of $\bC$ is equivalent to the category of modules of $D(H)$ the Drinfeld (or quantum) double of $H$. For the details of Drinfeld doubles, see \cite[Section 10.3]{Montgomery1}, \cite[Section 4.2.D]{CP1}, \cite[Section IX.4]{Kassel1}, or \cite[Section 7.14]{EGNO1}. The Drinfeld double $D(H)$ is isomorphic as a vector space to $(H^{\op}) ^* \otimes H$, and contains both $H$ and $(H^{\op})^*$ as Hopf subalgebras. Here if $H$ is a Hopf algebra with multiplication $\mu$, unit $\eta$, comultiplication $\Delta$, counit $\epsilon$, and antipode $S$, then $(H^{\op})^*$ is the Hopf algebra with multiplication $\Delta^*$, unit $\epsilon^*$, comultiplication $(\mu^{\op})^*$, counit $\eta^*$, and antipode $(S^{-1})^*$. 

The following result of Etingof--Ostrik will be important in extending the forgetful functor $\drin(\bC) \to \bC$ to the stable categories \cite{EO1}. 
 
\bpr{drinf-center}
If $\bC$ is a finite tensor category, then its Drinfeld center $\drin(\bC)$ is a finite tensor category, and the forgetful functor $F$ is exact and sends projective objects to projective objects.
\epr

The fact that $F$ preserves projectivity is a generalization of the classical Nichols--Zoeller theorem for Hopf algebras, which states that a finite-dimensional Hopf algebra is free as a module over any Hopf subalgebra \cite{NZ1}.

\section{Drinfeld Centers and Balmer Spectra}
\label{s:DCBS}

In this section, we prove general results relating the Balmer spectrum of $\st(\bC)$ to the Balmer spectrum of $\st(\drin(\bC))$, under the assumption that $\bC$ is an arbitrary finite tensor category.

\subsection{Construction of a continuous map between Balmer spectra}
\label{ss:mapf}

Recall the stable categories defined in Section \ref{ss:stab}. For the rest of this section, let $\bC$ be a finite tensor category, $\st(\bC)$ its stable category, $\drin(\bC)$ its Drinfeld center, $\st ( \drin ( \bC))$ the stable category of its Drinfeld center (which may be formed by \prref{drinf-center}), and $\St(\bC)$ and $\St(\drin(\bC))$ the respective ``big" stable categories, recall Notation \ref{big-st-name}. We have a forgetful functor $F: \drin( \bC) \to \bC$, and we have functors $G: \bC \to \st ( \bC)$ and $H: \drin(\bC) \to \st ( \drin ( \bC))$. The functor $F$ extends to a functor $\Ind(\drin(\bC)) \to \Ind(\bC)$, which we again denote by $F$, by \cite[Proposition 6.1.9]{KS1}.  We have the respective Balmer support data associated to $\st(\bC)$ and $\st (\drin(\bC))$: $$V_{\st \bC}: \st(\bC) \to \mc{X}_{\cl}(\Spc \st(\bC))$$ 
and $$V_{\st (\drin(\bC))}:\st( \drin( \bC)) \to \mc{X}_{\cl}(\Spc  \st (\drin (\bC))),$$ defined in their respective categories by sending $$A \mapsto \{\text{primes not containing }A\}.$$ 
\begin{notation}
\label{supp-z-c}
For readability, when $\bC$ is a finite tensor category we will denote $V_{\bC} := V_{\st \bC}$ and $V_{\drin}:=V_{\st (\drin(\bC))}$. The corresponding maps $\Phi$ (recalling the construction from Section \ref{ss:supp})  associated to these support data will similarly be denoted $\Phi_{\bC}$ and $\Phi_{\drin}$, respectively. We will similarly denote open sets in the Balmer spectrum on these respective categories by $U_{\bC}:=U_{\st(\bC)}$ and $U_{\drin}:=U_{\st(\drin(\bC))}$, recall the notation from Section \ref{ss:ttg}.
\end{notation}

The following proposition is probably well-known to experts, but we record it for completeness.

\bpr{funct-drinf-cent-st}
There is a functor $\overline{F} : \St(\drin(\bC)) \to \St (\bC)$ which extends the forgetful functor $F$, i.e. the diagram of functors
\begin{center}
\begin{tikzcd}
\Ind(\bC) \arrow[r, "G"]                     & \St(\bC)                                        \\
\Ind(\drin(\bC)) \arrow[u, "F"] \arrow[r, "H"] & \St(\drin(\bC)) \arrow[u, "\overline{F}", dashed]
\end{tikzcd}
\end{center}
commutes. This functor $\overline{F}$ is monoidal and triangulated.
\epr

\begin{proof}

Since the objects of $\St(\drin(\bC))$ are the in bijection with those of $\Ind(\drin(\bC))$, the functor $\overline{F}$ is well-defined on objects, namely by defining $$\overline{F} ( H(X)) := G(F(X)).$$ Let $f: X \to Y$ be a morphism in $\Ind(\drin(\bC))$. Then for $\overline{F}(H(f)) := G F (f)$ to be well-defined, we need $GF(g) =0$ for each morphism $g$ which factors through a projective in $\Ind(\drin(\bC))$. In other words, we need $F(g)$ to factor through a projective in $\Ind(\bC)$. Hence, to define $\overline{F}$, it is enough to know that $G \circ F$ sends all projective objects of $\Ind(\drin (\bC))$ to 0, which is true by \prref{drinf-center}.

Let $H(X) \in \St(\drin (\bC))$ an arbitrary object, where $X \in \Ind(\drin (\bC)).$ Then $\Sigma H(X)$ is defined as $H(Z)$, such that there exists a short exact sequence
$$0 \to X \to P \to Z \to 0$$ in $\Ind(\drin(\bC))$, where $P$ is a projective object in $\Ind(\drin(\bC))$. The object $\Sigma H(Z)$ is well-defined in $\St(\drin(\bC))$, by Schanuel's Lemma. Since $F$ is  exact and sends projectives to projectives, $$0 \to F(X) \to F(P) \to F(Z) \to 0$$ is an exact sequence in $\bC$ with $F(P)$ projective; therefore, $\Sigma ( G F( X)) \cong GF(Z)$ in $\st (\bC)$, and so we have $\overline{F}(\Sigma X) \cong \Sigma \overline{F}(X)$.

Now, let $X \to Y \to Z \to \Sigma X$ be a distinguished triangle in $\St(\drin(\bC))$. Then it is isomorphic to a triangle of the form
$$H(X') \to H(Y') \to H(Z') \to \Sigma H(X')$$ for some short exact sequence $$0 \to X' \to Y' \to Z' \to 0$$ in $\Ind(\drin ( \bC))$. Since $F$ is exact, and $G$ sends exact sequences to triangles, we have that the composition $GF$ is exact and hence
$$\overline{F}  H (X') \to \overline{F} H(Y') \to \overline{F} H(Z') \to \Sigma \overline{F} H(X')$$ is a triangle in $\St (\bC)$. Therefore, $$\overline{F} (X) \to \overline{F} (Y) \to \overline{F} (Z) \to \Sigma \overline{F} (X)$$ is a triangle as well, and so $\overline{F}$ is a triangulated functor. 
\end{proof}

For braided tensor triangulated categories, the Balmer spectrum $\Spc$ is functorial, as Balmer has shown in \cite[Proposition 3.6]{Balmer1}. This is a categorical reflection the ring-theoretic fact that $\Spec$ is functorial for commutative rings. On the other hand, for noncommutative rings, $\Spec$ is not a functor (for an in-depth exploration of the extent of the failure of functoriality of $\Spec$ for noncommutative rings, see \cite{Reyes1}). It is not suprising, then, that for generic monoidal triangulated categories, the Balmer spectrum is also not functorial; in other words, an monoidal triangulated functor between monoidal triangulated categories does not necessarily induce a map between their Balmer spectra. 

However, reflecting the classical prime ideal contraction for noncommutative rings, the forgetful functor $\overline{F}$ does induce a map between the Balmer spectra of $\st (\bC)$ and $\st(\drin(\bC))$. 

\bpr{fbar-induces-spc}
The functor $\overline{F}$ induces a continuous map $\Spc \st(\bC) \xrightarrow{f} \Spc  \st(\drin(\bC))$, defined by
$$f(\bP):= \{ X \in \st (\drin (\bC)): \overline{F}(X) \in \bP\}.$$
\epr

\begin{proof}
Let $\bP$ be a prime ideal of $\st(\bC)$. We must first show that $f(\bP)$ is a prime ideal of $\st(\drin(\bC))$.  

We first check that $f(\bP)$ is a thick ideal of $\st(\drin(\bC))$. This necessitates checking three properties:

(Triangulated) Suppose $\Sigma X \in f(\bP)$, in other words, $\overline{F}(\Sigma X) \in \bP$. Since $\overline{F}$ is triangulated, this is true if and only if $\Sigma \overline{F} (X) \in \bP$, which is true if and only if $\overline{F}(X) \in \bP$, in other words, $X \in f(\bP)$. Now, suppose $$X \to Y \to Z \to \Sigma X$$ is a distinguished triangle with $X$ and $Y$ in $f(\bP)$. This means that $\overline{F} (X)$ and $\overline{F}(Y)$ are in $\bP$. Since $\overline{F}$ is triangulated, the triangle
$$\overline{F}(X) \to \overline{F}(Y) \to \overline{F}(Z) \to \Sigma \overline{F}(X)$$ is distinguished in $\st(\bC)$. Now since the first two objects are in $\bP$, so is $\overline{F}(Z)$, and so $Z \in f(\bP)$.

(Thick) If $X \oplus Y$ is in $f(\bP)$, then $\overline{F}(X \oplus Y) \in \bP$; $\overline{F}$ is an additive functor, and so $\overline{F}(X) \oplus \overline{F}(Y) \in \bP$. This implies that both $\overline{F}(X)$ and $\overline{F}(Y)$ are in $\bP$, and so $X$ and $Y$ are both in $f(\bP)$.

(Ideal) Suppose $X \in f(\bP)$ and $Y \in \st(\drin(\bC))$. Since $\overline{F}$ is monoidal, we have $\overline{F}(X \otimes Y) \cong \overline{F}(X) \otimes \overline{F}(Y)$. Since $\overline{F}(X) \in \bP$, so is $\overline{F}(X) \otimes \overline{F}(Y)$, and thus $\overline{F}(X \otimes Y) \in \bP$ as well. Hence $X \otimes Y \in f(\bP)$. The symmetric argument shows that $Y \otimes X$ is in $f(\bP)$ as well, so $f(\bP)$ is a two-sided ideal.

(Prime) Let $A \otimes B \in f(\bP)$. Then $\overline{F} (A) \otimes \overline{F}(B) \in \bP$. But $\overline{F}(A)$ and $\overline{F}(B)$ commute with every object of $\st (\bC)$: by the ideal property of $\bP$, we have 
\begin{align*}
&\st(\bC) \otimes \overline{F}(A) \otimes \overline{F}(B) \subseteq \bP\\
 &\Rightarrow \overline{F}(A) \otimes \st(\bC) \otimes \overline{F}(B) \subseteq \bP\\
 &\Rightarrow \overline{F} (A) \text{ or } \overline{F}(B) \in \bP,
\end{align*}
 with the last step following by primeness of $\bP$. This implies that either $A$ or $B$ is in $f(\bP)$, which means that $f(\bP)$ is prime. 

We can also check directly that $f$ is continuous: an arbitrary closed set of $\Spc ( \st (\drin(\bC)))$ is of the form $V_{\drin}(\mc{T})= \{ \bP \in \Spc ( \st (\drin(\bC))) : \mc{T} \cap \bP = \varnothing\}$ (recalling Notation \ref{supp-z-c}) for some collection of objects $\mc{T}$ of $\st(\drin(\bC))$. Then
\begin{align*}
f^{-1} ( V_{\drin}(\mc{T})) &= \{ \bP \in \Spc \st(\bC): \mc{T} \cap \{ X \in \st (\drin (\bC)): \overline{F}(X) \in \bP\} = \varnothing\}\\
&= \{ \bP \in \Spc \st(\bC) : \overline{F} (\mc{T}) \cap \bP = \varnothing \}\\
&= V_{\bC} ( \overline{F} (\mc{T})),
\end{align*}
where by $\overline{F}(\mc{T})$ we mean the collection $\{\overline{F}(X) : X \in \mc{T} \}$.
\end{proof}

\bre{hopf-intro}
Recall the construction of the Drinfeld double from Section \ref{ss:drinf}. If $R$ is a finite-dimensional Hopf algebra, then $\drin(\modd(R)) \cong \modd( D(R))$. In this case, $D(R) \cong D((R^{\op})^*)$, and so we have two functors:
\begin{center}
\begin{tikzcd}
\modd(R) &                                                                                                    & \modd((R^{\op})^*) \\
         & \drin(\modd(R)) \cong \modd(D(R)) \cong \drin(\modd((R^{\op})^*)) \arrow[lu, "F_R"] \arrow[ru, "F_{(R^{\op})^*}"'] &           
\end{tikzcd}
\end{center}
which then give two maps between Balmer spectra:
\begin{center}
\begin{tikzcd}
\Spc(\stmod(R)) \arrow[rd, "f_R"'] &                   & \Spc(\stmod((R^{\op})^*)) \arrow[ld, "f_{(R^{\op})^*}"] \\
                                  & \Spc(\stmod(D(R))) &                                       
\end{tikzcd}
\end{center}
\ere

\subsection{A support data interpretation}
\label{ss:support}

We can interpret the map $f$ in the context of support data (recalling the definition from Section \ref{ss:supp}), by first defining a new support datum given as the composition of the functor $\overline{F}$ with the Balmer support $V_{\bC}$ on $\st(\bC)$. 

\bpr{support-datum}
Define a map $W: \st (\drin ( \bC)) \to \mc{X}_{\cl}(\Spc \st (\bC))$ by 
$$W ( X):= V_{ \bC} (\overline{F}(X))=\{ \bP \in \Spc \st (\bC) : \overline{F} (X) \not \in \bP\}.$$
This map is a support datum.
\epr

\begin{proof}
The first four conditions follow directly from the facts that $\overline{F}$ is a triangulated functor and $V_{\bC}$ is itself a support datum, since
\begin{enumerate}
\item $\overline{F} (0_{\st(\drin(\bC))} ) = 0_{\st(\bC)}$,
\item $\overline{F} ( X \oplus Y) = \overline{F} (X) \oplus \overline{F}(Y)$,
\item $\overline{F} ( \Sigma X) \cong \Sigma \overline{F}(X)$,
\item and if $X \to Y \to Z \to \Sigma X$ is a distinguished triangle, then so is $\overline{F}(X) \to \overline{F}(Y) \to \overline{F}(Z) \to \Sigma \overline{F}(X)$.
\end{enumerate} 

To check the last condition, we need to show that $$\bigcup_{Z \in \st ( \drin (\bC))} W ( X \otimes Z \otimes Y) = W(X) \cap W(Y).$$
By the ideal condition, if $\bP$ is a prime ideal which does not contain $\overline{F}(X) \otimes \overline{F}(Z) \otimes \overline{F}(Y)$ for some object $Z$, then it must also not contain $\overline{F}(X)$ or $\overline{F}(Y)$. Hence, $$\bigcup_{Z \in \st ( \drin (\bC))} W ( X \otimes Z \otimes Y) \subseteq W(X) \cap W(Y)$$ is automatic. 

For the reverse containment, suppose $\bP$ is a prime ideal which does not contain $\overline{F}(X)$ or $\overline{F}(Y)$. By the prime condition, that means $\bP$ does not contain the entire collection of objects $\overline{F}(X) \otimes \st(\bC) \otimes \overline{F}(Y)$. But since $\overline{F}(X)$ and $\overline{F}(Y)$ commute up to isomorphism with all elements of $\st (\bC)$, if $\overline{F}(X) \otimes \overline{F}(Y) \in \bP$, that would imply there is a containment $\overline{F}(X) \otimes \overline{F}(Y) \otimes \st(\bC) \subseteq \bP$, which would then imply $\overline{F}(X) \otimes \st(\bC) \otimes \overline{F}(Y) \subseteq \bP$, a contradiction. Hence, $\bP \in W ( X \otimes Y)$, and we have the claimed equality.
\end{proof}

By the universal property of the Balmer spectrum as in \thref{univ-balm}, the support datum $W$ induces a continuous map $\Spc \st (\bC) \to \Spc \st (\drin (\bC))$. This map is defined as
$$\bP \mapsto \{ X \in \st (\drin (\bC)) : \bP \not \in W(X)\}$$ by \cite[Theorem 4.2.2]{NVY1}. One may observe that this map is the same as the map defined in \prref{fbar-induces-spc}. We have the following diagram, which commutes by definition:
\begin{center}
\begin{tikzcd}
\st(\drin(\bC)) \arrow[rd, "\overline{F}"'] \arrow[rr, "W"] &                                 & \mc{X}_{\cl} ( \Spc ( \st( \bC))) \\
                                                            & \st(\bC) \arrow[ru, "V_{\bC}"'] &                                 
\end{tikzcd}
\end{center}
On the level of ideals, we have the following induced maps, recall $\Phi$ and $\Theta$ associated to a support datum as constructed in Section \ref{ss:supp}:
\begin{center}
\begin{tikzcd}
\ThickId(\st(\drin(\bC))) \arrow[rr, "\Phi_W"] \arrow[rd, "\Psi"'] &                                                                               & \mc{X}_{\Sp}(\Spc(\st(\bC))) \arrow[ll, "\Theta_W"', bend right] \arrow[ld, "\Theta_{\bC}", bend left] \\
                                                                   & \ThickId(\st(\bC)) \arrow[ru, "\Phi_{\bC}"'] \arrow[lu, "\Lambda", bend left] &                                                                                                      
\end{tikzcd}
\end{center}
Here, for thick ideals $\bI$ of $\st(\drin(\bC))$ and $\bJ$ of $\st(\bC)$, the maps $\Psi$ and $\Lambda$ are defined by
\begin{align*}
\Psi: \bI &\mapsto \langle \overline{F}(\bI) \rangle,\\
\Lambda: \bJ & \mapsto \{ X \in \st(\drin(\bC)): \overline{F}(X) \in \bJ \}.
\end{align*}
By definition, the inner and outer triangles commute: in other words, $\Phi_W= \Phi_{\bC} \circ \Psi$, and $\Theta_W = \Lambda \circ \Theta_{\bC}$. 

\subsection{Recovering ideals from their supports}
\label{ss:mapf}

In \cite[Theorem 6.2.1]{NVY1}, conditions were given under which an arbitrary support datum $\sigma: \bT \to \mc{X}(S)$ has the property that $\Phi_{\sigma}$ is a left, right, and two-sided inverse to $\Theta_{\sigma}$. If $\Phi_{\sigma}$ is a left inverse to $\Theta_{\sigma}$, this means that all thick ideals can be recovered from their supports; when $\Phi_{\sigma}$ and $\Theta_{\sigma}$ are a mutually inverse bijection, the ideals are completely classified by the topological space $S$. Since the support datum $W(-)$ defined above might not satisfy conditions under which every ideal may be recovered from their support (see Section \ref{s:applications} for examples), in this section we discuss precisely which ideals can be recovered in this way; this allows us to describe the image of the map $f$ defined above.

We now introduce some terminology, which will be useful for our reconstruction theory.

\begin{notation}
\label{kernel}
When the finite tensor category $\bC$ is clear by context, we will denote by $\bK$ the kernel of the functor $\overline{F}: \St(\drin(\bC)) \to \St(\bC)$. 
\end{notation}

An equivalent characterization of the kernel of $\overline{F}$ can be given by
$$\bK = \{H(X) : X \in \Ind(\drin(\bC)) \text{ such that }F(X) \text{ is projective in }\Ind(\bC)\}.$$
This follows from the fact that the objects of $\St(\bC)$ isomorphic to 0 correspond precisely to the projective objects of $\Ind(\bC)$, as we saw in Section \ref{ss:stab}. 

\ble{hk-ideal}
The kernel of $\overline{F}$ is a thick localizing ideal of $\St(\drin(\bC))$.
\ele

\begin{proof}
Since $\overline{F}$ is a monoidal triangulated functor, it is straightforward to verify that the collection of objects $X$ such that $\overline{F}(X) \cong 0$ is closed under taking cones, shifts, direct summands, and by tensoring on the left or right by arbitrary objects of $\St(\drin(\bC))$. The functor $F$ commutes with arbitrary coproducts by \cite[Proposition 6.1.9]{KS1}, and so the kernel of $\overline{F}$ is closed under arbitrary coproducts, i.e. $\bK$ is localizing.  
\end{proof}

\ble{emptyset-kernel}
An object $A \in \st(\drin(\bC))$ satisfies $W(A) = \varnothing$ if and only if $A \in \bK$. 
\ele

\begin{proof}
First, note that if $A \in \bK$, then by definition $\overline{F}(A) \cong 0$, and so $$W(A) = V_{\bC} (0) =\{ \bP \in \Spc ( \st(\bC)) : 0 \not \in \bP\} = \varnothing.$$ 

For the other direction, recall that by the rigidity of $\bC$, all thick ideals of $\st (\bC)$ are semiprime, i.e. intersections of prime ideals, by \prref{semiprime}. This implies in particular that the ideal $\langle 0 \rangle$ is semiprime; in other words, the only object contained in all prime ideals of $\st(\bC)$ is 0. By definition, this means that if $X$ is an object of $\st (\bC)$ such that $V_{\bC} (X) = \varnothing$, then $X\cong 0$. Hence, we have
$$\varnothing= W (A)= V_{\bC} ( \overline{F}(A)) \Rightarrow \overline{F}(A)  \cong 0 \Rightarrow A \in \bK.$$
\end{proof}

Using the localization and colocalization functors defined in Section \ref{ss:compact-gen}, we are now able to prove the following, which is the critical step in determining which ideals can be recovered from their $W$-support and determining the image of the map $f: \Spc \st(\bC) \to \Spc \st(\drin(\bC))$ defined in \prref{fbar-induces-spc}.

\bth{ideal-recover}
Let $\bI$ be a thick ideal of $\st(\drin(\bC))$ such that $\Loc(\bI)$ contains $\bK$. Suppose that $X$ is an object of $\st(\drin(\bC))$ such that $\overline{F}(X) \in \langle \overline{F} (\bI) \rangle$, that is, the thick ideal of $\st(\bC)$ generated by all $\overline{F}(Y)$ for $Y \in \bI$. Then $X$ is in $\bI$.
\eth

\begin{proof}
By \thref{local}, we have a distinguished triangle 
$${\it{\Gamma}}_{\bI} (X) \to X \to L_{\bI} (X) \to \Sigma {\it{\Gamma}}_{\bI} (X)$$
in $\St(\drin(\bC))$, using the localization and colocalization functors associated to the thick ideal $\bI$. We know that there are no morphisms from $\bI$ to $L_{\bI} (X)$; in other words, if $Y \in \bI$ and $Z$ is any compact object in $\St(\drin(\bC))$, then 
\begin{align*}
0&=\Hom_{\St(\drin(\bC))}(Z \otimes Y , L_{\bI} (X)) \\
&\cong \Hom_{\St(\drin(\bC))}(Z, L_{\bI} (X) \otimes Y^*).
\end{align*}
Since this holds for all compact objects $Z$, this implies that $L_{\bI} (X) \otimes Y^* \cong 0$. Since all compact objects are rigid, and by \prref{semiprime} all thick ideals are closed under taking duals, we have $L_{\bI} (X) \otimes Y \cong 0$ for all $Y \in \bI$. Since $\overline{F}$ is a monoidal functor, this additionally implies that $$\overline{F} (L_{\bI} (X)) \otimes \overline{F} (Y) \cong 0$$ in $\St(\bC)$, for all $Y \in \bI$. 

Now, consider the thick ideal $\langle \overline{F} (\bI) \rangle$. This is formed successively by taking shifts, cones, direct summands, and tensor products with arbitrary elements of $\st(\bC)$, starting from the collection of objects of the form $\overline{F}(Y)$ for $Y \in \bI$. This allows us to conclude inductively that $\overline{F} ( L_{\bI} (X)) \otimes A \cong 0$ for all $A$ in $\langle \overline{F} (\bI) \rangle$, since inductively each step by which we construct $\langle \overline{F} (\bI) \rangle$ preserves the property that tensoring with $\overline{F} (L_{\bI} (X))$ gives 0. To be more explicit, if
$$A \to B \to C \to \Sigma A$$ is a distinguished triangle in $\st(\bC)$ such that $A \otimes  \overline{F} (L_{\bI} (X)) \cong B \otimes \overline{F} (L_{\bI} (X)) \cong 0$, then it is straightforward that additionally $C \otimes \overline{F} (L_{\bI} (X)) \cong 0$ as well. Similarly, if $A \otimes \overline{F} (L_{\bI} (X)) \cong 0$, then $\Sigma (A) \otimes \overline{F} (L_{\bI} (X)) \cong \Sigma ( A \otimes \overline{F} (L_{\bI} (X))) \cong \Sigma 0 \cong 0$. Furthermore, if we have $(A \oplus B) \otimes \overline{F} (L_{\bI} (X)) \cong 0$, then we also have $A \otimes \overline{F} (L_{\bI} (X)) \cong 0 \cong B \otimes \overline{F} (L_{\bI} (X))$. Lastly, if $A \otimes \overline{F} (L_{\bI} (X)) \cong 0$ and $B$ is an arbitrary object in $\st(\bC)$, then $A \otimes B \otimes \overline{F} (L_{\bI} (X)) \cong A \otimes \overline{F} (L_{\bI} (X)) \otimes B \cong 0$ as well, using the commutativity of $\overline{F} (L_{\bI} (X))$.

To reiterate, the upshot of the previous paragraph is that $A \otimes \overline{F} (L_{\bI} (X)) \cong 0$ for all $A \in \langle \overline{F} (\bI) \rangle$. But by assumption, we have $\overline{F} (X) \in \langle \overline{F} (\bI) \rangle$. Hence, $$ \overline{F} ( X \otimes L_{\bI} (X)) \cong \overline{F} (X) \otimes \overline{F} (L_{\bI} (X)) \cong 0.$$
Therefore, $X \otimes L_{\bI} (X) $ is an object in $\bK$, the collection of objects of $\St(\drin(\bC))$ mapped to 0 by $\overline{F}$. By assumption, $\Loc(\bI)$ contains $\bK$, and so $X \otimes L_{\bI} (X) \in \Loc(\bI)$. 

Now, consider the distinguished triangle obtained by tensoring the triangle
$${\it{\Gamma}}_{\bI} (X) \to X \to L_{\bI} (X) \to \Sigma {\it{\Gamma}}_{\bI} (X)$$
by $X$: this gives us 
$$X \otimes {\it{\Gamma}}_{\bI} (X)  \to X \otimes X \to X\otimes  L_{\bI} (X) \to \Sigma X \otimes  {\it{\Gamma}}_{\bI} (X).$$
We have just finished showing that the third object of this triangle is in $\Loc(\bI)$. The first object is in $\Loc(\bI)$ as well, by \thref{local}. Since $\Loc(\bI)$ is triangulated, this implies $X \otimes X$ is in $\Loc(\bI)$. But by \cite[Lemma 2.2]{Neeman2}, since $\bI$ is a thick subcategory of compact objects, the compact objects in $\Loc(\bI)$ are precisely the objects of $\bI$. Thus, $X \otimes X \in \bI$, and by semiprimeness of $\bI$ (\prref{semiprime}) so is $X$; this completes the proof.
\end{proof}

We can now give a condition under which an ideal $\bI$ can be recovered from its support $\Phi_W(\bI)$.

\bco{theta-phi-inv}
Let $\bI$ be an ideal such that $\Loc(\bI)$ contains $\bK$. Then $\Theta_W \circ \Phi_W ( \bI) = \bI$.
\eco

\begin{proof}
By definition, 
\begin{align*}
\Theta_W \circ \Phi_W (\bI) &= \Theta_W ( \Phi_{\bC} ( \overline{F} (\bI))) \\
&= \{ X \in \st(\drin(\bC)) : W(X) \subseteq \Phi_{\bC} (\overline{F} (\bI))\}\\
&= \{ X \in \st(\drin(\bC)) : V_{\bC} ( \overline{F} (X)) \subseteq \Phi_{\bC} (\langle \overline{F} (\bI) \rangle )\}\\
&= \{ X \in \st(\drin(\bC)) : \; \forall \;\bP \in \Spc \st(\bC) \text{ with }\overline{F}(X) \not \in \bP, \langle \overline{F} (\bI) \rangle \not \subseteq \bP \}\\
&= \{ X \in \st(\drin(\bC)) : \; \forall \; \bP \in \Spc \st(\bC) \text{ with } \langle \overline{F} (\bI) \rangle \subseteq \bP, \overline{F}(X)\in \bP \}\\
&= \left \{ X \in \st(\drin(\bC)) : \overline{F} (X) \in \bigcap_{\bP \in \Spc \st(\bC), \langle \overline{F}(\bI) \rangle \subseteq \bP} \bP \right \}\\
&= \{X \in \st(\drin(\bC)): \overline{F}(X) \in \langle \overline{F}(\bI) \rangle \}.
\end{align*}
The last equality follows from \prref{semiprime}. The corollary now follows directly from \thref{ideal-recover}.
\end{proof}

\subsection{The image of prime ideal contraction}
\label{ss:imf}

We now describe the relationship of the image of the map $f$ to the kernel $\bK$ of $\overline{F}$. 

\bpr{map-f-surj}
Let $\bC$ be a finite tensor category.
\begin{enumerate}
\item If $\bP$ is in the image of the map $f: \Spc \st ( \bC) \to \Spc  \st ( \drin (\bC))$, then $\bP$ contains $\bK \cap \st(\drin(\bC))$, the kernel of $\overline{F}$ restricted to compact objects. 
\item If $\bP$ is a prime ideal of $\st(\drin(\bC))$ such that $\Loc(\bP)$ contains $\bK$, then $\bP$ is in the image of $f$.
\end{enumerate}
\epr

\begin{proof}
For (1), if $\bQ$ is a prime ideal of $\st(\bC)$, then $f(\bQ)$ contains $\bK \cap \st(\drin(\bC))$, which are by definition the finite-dimensional objects $X$ such that $\overline{F}(X) \cong 0$: if $X$ is in $\st(\drin(\bC))$ and $\overline{F}(X) \cong 0$, then $X \in \{Y \in \st(\drin(\bC)) : \overline{F}(Y) \in \bQ\} = f(\bQ)$, since 0 is in every prime ideal of $\st(\bC)$.

Part (2) is an application of both \thref{ideal-recover} and \cite[Theorem 3.2.3]{NVY1}. Let $\bP$ be a prime ideal of $\st ( \drin ( \bC))$ such that $\Loc(\bP)$ contains $\bK$. Consider the following two collections of objects in $\st(\bC)$:
\begin{enumerate}
\item the ideal $\bI := \langle \overline{F}(X) : X \in \bP \rangle$ of $\st(\bC);$
\item the collection $\mc{M} := \{ \overline{F} ( Y) : Y \not \in \bP \}$ of objects in $\st(\bC)$.
\end{enumerate}

We first claim that these two collections of objects are disjoint. If $\overline{F}(Y) \in \bI$ then $Y \in \Theta_W (\Phi_W ( \bP))$, implying that $Y \in \bP$ by \coref{theta-phi-inv}. This means that in particular, if $\overline{F}(X) \cong \overline{F}(Y)$, then either both $X$ and $Y$ are in $\bP$, or neither are, and so $\bI$ and $\mc{M}$ are indeed disjoint. 

Since $\bP$ is a proper ideal of $\st(\drin(\bC))$, it follows that $\mc{M}$ is nonempty, and thus $\bI$ is a proper ideal of $\st(\bC)$. We claim that $\mc{M}$ is a multiplicative subset. Suppose $\overline{F} (X)$ and $\overline{F} (Y)$ are in $\mc{M}$. Then if $\overline{F} (X) \otimes \overline{F} (Y) \cong \overline{F} (X \otimes Y)$ was not in $\mc{M}$, this would imply that $X \otimes Y \in \bP$; by the prime condition of $\bP$, either $X$ or $Y$ (without loss of generality, say $Y$) would then be in $\bP$. This is a contradiction, since $\overline{F} (Y) \in \mc{M}$ implies $Y \not \in \bP$, which is a consequence of the observation above that $\bI$ and $\mc{M}$ are disjoint.

By \cite[Theorem 3.2.3]{NVY1}, given a disjoint pair consisting of a multiplicative subset and a proper ideal of any monoidal triangulated category (in this case, $\st(\bC)$), there exists a prime ideal $\bQ$ of $\st(\bC)$ such that $\bQ \cap \mc{M} = \varnothing$ and $\bI \subseteq \bQ$. We have 
$$f(\bQ)= \{ X \in \st (\drin(\bC)) : \overline{F} (X) \in \bQ\},$$ and then since $\bI \subseteq \bQ$, it is automatic that $\bP \subseteq f(\bQ)$; and since $\bQ$ is disjoint from $\mc{M}$, in fact $\bP = f(\bQ)$. Thus, $f$ surjects onto the collection of prime ideals $\bP$ such that $\Loc(\bP)$ contains $\bK$, which completes the proof.
\end{proof}

By \prref{map-f-surj}, we have inclusions of the following subsets of $\Spc \st(\drin(\bC))$:
\begin{equation}
\label{ee:ineq-imf}
\{ \bP : \bK \cap \st(\drin(\bC)) \subseteq \bP \} \supseteq \im f \supseteq \{ \bP  : \bK \subseteq \Loc(\bP)\}.
\end{equation}

We note the following lemma, which is a special case of \cite[Proposition 1.47]{BIK2}.

\ble{k-comp-equiv}
The following are equivalent.
\begin{enumerate}
\item The kernel $\bK$ of $\overline{F}$ is generated as a localizing category (recalling Section \ref{ss:compact-gen}) by the set $\bK \cap \st(\drin(\bC))$.
\item For every nonzero $X$ in $\bK$, there exists a compact object $Y$ in $\bK$ which has some nonzero map $Y \to X$ in $\St(\drin(\bC))$.
\end{enumerate}
\ele

In particular, to prove that (2) $\Rightarrow$ (1), one simply observes that if $X \in \bK$, then the distinguished triangle
\[
{\it{\Gamma}}_{\bK \cap \st(\drin(\bC))} X \to X \to L_{\bK \cap \st(\drin(\bC))} X\to \Sigma {\it{\Gamma}}_{\bK \cap \st(\drin(\bC))} X
\]
given by \thref{local} implies that $L_{\bK \cap \st(\drin(\bC))} X \in \bK$. But by definition, it is in the perpendicular space to $\bK \cap \st(\drin(\bC))$, which by the assumption of (2) means that it is 0. Hence ${\it{\Gamma}}_{\bK \cap \st(\drin(\bC))} X \cong X$, that is, $X$ is in $\Loc(\bK \cap \st(\drin(\bC)))$.

If these conditions are satisfied, then we can sharpen (\ref{ee:ineq-imf}), as well as \coref{theta-phi-inv}.

\bco{map-f-kernel-cg}
Suppose the kernel $\bK$ of $\overline{F}$ satisfies the equivalent conditions of \leref{k-comp-equiv}. 
\begin{enumerate}
\item The image of $f$ is precisely the collection of prime ideals of $\st(\drin(\bC))$ which contain $\bK \cap \st(\drin(\bC))$, that is, the collection of objects $X$ in $\st(\bC)$ such that $\overline{F}(X) \cong 0$.
\item A thick ideal $\bI$ of $\st (\drin(\bC))$ satisfies $\Theta_W \circ \Phi_W(\bI) = \bI$ if and only if $\bI$ contains $\bK \cap \st(\drin(\bC))$. 
\end{enumerate}
\eco

\begin{proof}
Suppose $\Loc( \bK \cap \st(\drin(\bC))) = \bK$. For (1), let $\bP$ be a prime ideal of $\st(\drin(\bC))$ containing $\bK\cap \st(\drin(\bC))$. Then $\Loc(\bP)$ contains $\Loc(\bK \cap \st(\drin(\bC))) = \bK$. Hence the collection of inequalities of (\ref{ee:ineq-imf}) becomes an equality, and we are done.

For (2), similarly, we have by \coref{theta-phi-inv} that if $\Loc(\bI)$ contains $\bK$, then $\Theta_W \circ \Phi_W(\bI) = \bI$. Since $\bK = \Loc(\bK \cap \st(\drin(\bC)))$, we have $\bK \subseteq \Loc(\bI)$ if and only if there is containment $\bK \cap \st(\drin(\bC)) \subseteq \bI$. For the other direction, we note that for any ideal $\bI$, we have $\bK \cap \st(\drin(\bC)) \subseteq \Theta_W \circ \Phi_W (\bI)$, and so any thick ideal satisfying $\Theta_W \circ \Phi_W (\bI)= \bI$ must have $\bK \cap \st(\drin(\bC)) \subseteq \bI$ as well. 
\end{proof}

\bre{sp-closed}
\coref{map-f-kernel-cg} implies that if $\bC$ satisfies the conditions of \leref{k-comp-equiv}, then the image of $f: \Spc  \st(\bC) \to \Spc \st(\drin(\bC))$ is automatically the complement of a specialization-closed set, since we have
\begin{align*}
\im(f) &= \{ \bP \in \Spc  \st(\drin(\bC)) : \bP \supseteq \bK \cap \st(\drin(\bC)) \} \\
&= \Spc \st(\drin(\bC)) \backslash (\Phi_{\drin} ( \bK \cap \st(\drin(\bC)) ) ).
\end{align*}
In other words, the image of $f$ can be written as an intersection of open sets. If $\bK~\cap~\st(\drin(\bC))$ is generated (as a thick ideal) by a finite collection of objects, say $\{X_i\}_{i=1}^n$, then it follows that $\im(f)$ is in fact an open subset of $\Spc \st(\drin(\bC))$, namely $$\im(f) =   U_{\drin}( X_1 \oplus... \oplus X_n )$$
(recall the notation of $U_{\drin}$ from Section \ref{ss:ttg} and Notation \ref{supp-z-c}).
\ere

\bre{ineq-not-eq}
In the situation of \coref{map-f-kernel-cg} (2), we have \coref{theta-phi-inv} sharpened from a one-way implication to a two-way implication. We note on the other hand that if the conditions of \leref{k-comp-equiv} are not satisfied, then \coref{theta-phi-inv} can never be an if-and-only-if, for the following reason. The collection of objects $\bK \cap \st(\drin(\bC))$ is itself a thick ideal of $\st(\drin(\bC))$, since it is in particular the kernel of the monoidal triangulated functor $\overline{F}$ restricted to compact objects. But now note that
\begin{align*}
\Theta_W \circ \Phi_W ( \bK \cap \st(\drin(\bC)) ) &= \{ X \in \st (\drin(\bC)): W(X) \subseteq \Phi_W ( \bK \cap \st(\drin(\bC)) ) \}\\
&= \{ X \in \st (\drin(\bC)): W(X) \subseteq \varnothing \}\\
&= \{ X \in \st(\drin(\bC)) : \overline{F} (X) \cong 0 \}\\
&= \bK \cap \st(\drin(\bC)).
\end{align*}
Here the first equality is by the definition of $\Theta_W$, the second and third equalities are by \leref{emptyset-kernel} and the definition of $\Phi_W$, and the last equality by the definition of the kernel $\bK$. In other words, the thick ideal $\bK \cap \st(\drin(\bC))$ can be recovered from its support. But plainly, since we are assuming the conditions of \leref{k-comp-equiv} are not satisfied, we have
$$\Loc(\bK \cap \st(\drin(\bC))) \not \supseteq \bK,$$
and so \coref{theta-phi-inv} cannot be sharpened to an if-and-only-if statement. 
\ere

\subsection{Conditions under which $f$ is injective, surjective, or a homeomorphism}
\label{ss:fhomeo}

We now give conditions under which $\Phi_W$ and $\Theta_W$ are inverses, and $f$ is surjective, injective, and a homeomorphism.

\bth{f-inj-surj-homeo}
Let $\bC$ be a finite tensor category.
\begin{enumerate}
\item  The following conditions are equivalent.
	\begin{enumerate}
	\item For all $X \in \bK$, there exists an isomorphism $X\cong 0$ in $\St(\drin(\bC))$.
	\item The map $f$ is surjective and $\bK$ is generated as a localizing category by its subcategory $\bK \cap \st(\drin(\bC))$.
	\item As maps $\ThickId(\st(\drin(\bC))) \to \ThickId(\st(\drin(\bC)))$, we have $\Lambda \circ \Psi= \id$. 
	\item As maps $\ThickId(\st(\drin(\bC))) \to \ThickId(\st(\drin(\bC)))$, we have $\Theta_W \circ \Phi_W= \id$.
	\end{enumerate}
\item If $\bC$ is braided, then the following hold.
	\begin{enumerate}
	\item The map $f$ is injective.
	\item As maps $\ThickId(\st(\bC)) \to \ThickId(\st(\bC))$, we have $\Psi \circ \Lambda = \id$. 
	\item If additionally $\Spc \st(\bC)$ is topologically Noetherian, then $\Phi_W \circ \Theta_W = \id$.
	\end{enumerate}
\item If $X\cong 0$ in $\St(\drin(\bC))$ for all $X \in \bK$ and $\bC$ is braided, then the following hold.
	\begin{enumerate}
	\item The map $f$ is a homeomorphism.
	\item The maps $\Psi$ and $\Lambda$ define mutually inverse bijections between $\ThickId(\st(\drin(\bC)))$ and $\ThickId(\st(\bC))$. 
	\item If additionally $\Spc \st(\bC)$ is topologically Noetherian, then $\Phi_W$ and $\Theta_W$ are mutually inverse bijections between $\ThickId(\st(\drin(\bC)))$ and $\mc{X}_{\Sp}(\Spc(\st(\bC)))$.
	\end{enumerate}
\end{enumerate}
\eth

\begin{proof}
Suppose (1a) holds, and so $\bK$ consists only of objects isomorphic to 0, in other words, for all objects $X \in \drin(\bC)$,
$$F(X) \text{ is projective in } \bC \Leftrightarrow X \text{ is projective in }\drin(\bC).$$ In particular this means that $\bK$ is generated by $\bK \cap \st(\drin(\bC))$, since all objects of $\bK$ are isomorphic to 0. Then (1c) follows from \thref{ideal-recover}, and the conditions (1b) and (1d) follow directly from \coref{map-f-kernel-cg}.

Now, suppose (1b) is satisfied. By \prref{map-f-surj}, this means that every prime ideal of $\st(\drin(\bC))$ contains $\bK \cap \st(\drin(\bC))$. But since every ideal is semiprime, the zero ideal is equal to the intersection of all primes of $\st(\drin(\bC))$, and so $\bK \cap \st(\drin(\bC))$ is contained in the zero ideal. Since $\bK$ is generated by $\bK \cap \st(\drin(\bC))$, i.e. the zero ideal, this implies that (1a) holds.  

Note that
\begin{align*}
\Lambda ( \Psi ( \langle 0_{\st(\drin(\bC))} \rangle))&= \Lambda ( \langle 0_{\st(\bC)} \rangle)\\
&= \{X \in \st(\drin(\bC)): \overline{F}(X) \in \langle 0_{\st(\bC)} \rangle \}\\
&= \{X \in \st(\drin(\bC)): \overline{F}(X) \cong 0\}\\
&= \bK.
\end{align*}
Hence (1c) implies (1a). 

Lastly, suppose condition (1d) holds. This implies by \coref{theta-phi-inv} that $\bK \subseteq \Loc(\bI)$ for every thick ideal $\bI$; in particular, this means that $\bK$ is contained in the localizing category generated by 0, which consists only of objects isomorphic to 0. Hence, (1a) holds. 

To show (2), first note that if $\bC$ is braided with a braiding $\gamma$, then $\overline{F}$ is essentially surjective, since for any object $X$ in $\bC$, the pair $(X, \gamma_X)$ is an object of $\drin(\bC)$ and $\overline{F}$ sends $H( X, \gamma_X)$ to $G(X)$. Now, we note that if $\bP$ and $\bQ$ are prime ideals of $\st(\bC)$, then:
\begin{align*}
f(\bP) &= f(\bQ)\\
& \Updownarrow \\ 
\{ X \in \st(\drin(\bC)) : \overline{F} (X) \in \bP \} &= \{X \in \st(\drin(\bC)): \overline{F}(X) \in \bQ \} \\
& \Updownarrow \\
\forall \; X \in \st(\drin(\bC)), \; \overline{F} (X) \in \bP & \Leftrightarrow \overline{F} (X) \in \bQ \\
& \Updownarrow \\
\forall \; Y \in \st(\bC), \; Y \in \bP & \Leftrightarrow Y \in \bQ\\
& \Updownarrow\\
\bP &= \bQ.
\end{align*}
Hence, if $\bC$ is braided then (2a) follows.

Condition (2b) also follows directly from the fact that $\overline{F}$ is essentially surjective. 

For (2c), recall that by \thref{balmer-top}, $\Spc(\st(\bC))$ is Noetherian if and only if every closed set is of the form $V_{\bC} ( A)$ for some object $A \in \st(\bC)$. If $S$ is a specialization-closed set in $\Spc(\st(\bC))$, then by definition
\begin{align*}
\Phi_W ( \Theta_W ( S)) &= \Phi_W ( \{ X \in \st(\drin(\bC)) : W(X) \subseteq S\})\\
&= \bigcup_{X \in \Theta_W(S)} W(X)\\
& \subseteq S.
\end{align*}
For the other direction, we can write $S$ as a union of closed sets, say $S= \bigcup_{i \in I} S_i$, and by the Noetherianity of $\Spc (\st(\bC))$, there exist objects $A_i$ of $\st(\bC)$ such that $S_i = V_{\bC} (A_i)$. Since $\overline{F}$ is essentially surjective, we can pick $X_i \in \st(\drin(\bC))$ with $\overline{F}(X_i) = A_i$. Since $$W(X_i) = V_{\bC} ( A_i)= S_i \subseteq S,$$ we have by definition each $X_i$ is in $\Theta_W ( S)$. Therefore,
\begin{align*}
\Phi_W ( \Theta_W (S)) & \supseteq \bigcup_{i \in I} W ( X_i)\\
&= \bigcup_{i \in I} S_i\\
&= S.
\end{align*}
Thus $S = \Phi_W ( \Theta_W (S))$. 

Suppose the assumptions of (3). Then (3b) and (3c) follow immediately from parts (1) and (2). To show (3a), it is enough to show that $f$ is a closed map, by (1a) and (2a). Take an arbitrary closed set $V_{\bC}(\mc{T})$ in $\Spc \st(\bC)$. We claim that the image of $V_{\bC} (\mc{T})$ under $f$ is precisely $V_{\drin} (\hat{\mc{T}})$, where $\hat{\mc{T}} = \{ X \in \st(\drin(\bC)) : \overline{F} (X) \in \mc{T} \}.$ 

For the first direction, suppose $\bP \in V_{\bC}(\mc{T})$, in other words, $\bP \cap \mc{T} = \varnothing$. Since $f(\bP) = \{ X : \overline{F} (X) \in \bP\}$, this implies that for all $X \in f(\bP)$, we have $X \not \in \hat{\mc{T}}$. Therefore $f(\bP ) \cap \hat{\mc{T}} = \varnothing$, and so $f(\bP) \in V_{\drin} (\hat{\mc{T}})$. This shows $ f( V_{\bC} (\mc{T})) \subseteq V_{\drin} (\hat{\mc{T}})$. 

For the other containment, suppose $\bQ$ is a prime ideal of $\st(\drin(\bC))$ in $V_{\drin}(\hat{\mc{T}})$. Then $\overline{F} (X) \not \in \mc{T}$ for all $X \in \bQ$. Since $f$ is surjective, we can pick $\bP \in \Spc \st(\bC)$ with $f(\bP) = \bQ$, and for all $\overline{F} (X) \in \bP$, we must have $\overline{F}(X) \not \in \mc{T}$. Since $\overline{F}$ is essentially surjective, this implies $A \not \in \mc{T}$ for all $A \in \bP$, and so $\bP \cap \mc{T} = \varnothing,$ i.e. $\bP \in V_{\bC}(\mc{T}).$ This shows the other containment $ f( V_{\bC} (\mc{T})) \supseteq V_{\drin} (\hat{\mc{T}})$, and so we have equality.

Hence, $f$ sends the closed set $V_{\bC}(\mc{T})$ to the closed set $V_{\drin} (\hat{\mc{T}})$, and so it is a continuous, bijective, closed map, and therefore a homeomorphism.
\end{proof}

\section{Applications}
\label{s:applications}

The time has come for concrete applications of our theory.

\subsection{Group algebras and dual group algebras}
\label{ss:groupalg}

Let $G$ be a finite group, $\kk$ be an algebraically closed field of characteristic $p$ which divides the order of $G$, and $\kk G$ the group algebra of $G$ over $\kk$. Let $\bC=\modd(\kk G)$, a finite tensor category. The Drinfeld double $D(\kk G)$ is a Hopf algebra containing $\kk G$ and $(\kk G ^{\op})^*$ as Hopf subalgbras. We will denote the dual of the group algebra by $\kk [G]$, and in that case we can write $(\kk G^{\op})^* = \kk[G]^{\cop}.$ The collection
$$\{ p_g h : g, h \in G\}$$ is a $\kk$-basis of $D(\kk G)$, where the elements $\{p_g : g \in G\}$ refer to the basis of $\kk[G]^{\cop}$ dual to the standard basis of $\kk G$. The multiplication is determined by the relations
$$h p_g = p_{h g h^{-1}} h,$$ 
see for instance \cite[Section IX.4.3]{Kassel1}. 

\ble{group-projectives}
Let $G$ and $\kk$ be as above and $F: \Mod(D(\kk G)) \to \Mod(\kk G)$ the forgetful functor. Then if $F(P)$ is projective as a $\kk G$-module, then $P$ is projective as a $D(\kk G)$-module. 
\ele

\begin{proof}
A module for $D(\kk G)$ is a $\kk G$ module $M$ which is also a $G$-graded vector space, such that if $m \in M$ is a homogeneous element of degree $g$, then $h.m$ is homogeneous of degree $hgh^{-1}$. Suppose we have a short exact sequence $$0 \to A \to B \xrightarrow{t} C \to 0$$ of $D(\kk G)$-modules such that
$$0 \to F(A) \to F(B) \to F(C) \to 0$$ is a split short exact sequence of $G$-modules. We claim that the original sequence splits as $D(\kk G)$-modules. Pick a homogeneous basis $\{c_i\}$ of $C$ under the $G$-grading, where $c_i$ has degree $g_i$. Now pick a splitting $s: C \to B$. Define 
$\hat s ( c_i)= p_{g_i} s(c_i).$ This map is homogeneous with respect to the $G$-grading, and it is still a $G$-module map:
\begin{align*}
g \hat s ( c_i) &= g p_{g_i}  s(c_i)\\
&= p_{g g_i g^{-1}} g s(c_i)\\
&= p_{g g_i g^{-1}} s ( g c_i)\\
&= \hat s ( g c_i).
\end{align*}
Since on the basis $\{c_i\}$ we have 
\begin{align*}
t \circ \hat s(c_i) &=t ( p_{g_i}. s( c_i))\\
&= p_{g_i}. t s(c_i)\\
&= p_{g_i}. c_i\\
&= c_i, 
\end{align*}
we have that $\hat s$ is a splitting of $D(\kk G)$-modules. 

Now, to prove the original claim, suppose $F(P)$ is projective as a $G$-module. Since $F$ is exact, this means that for every short exact sequence
$$0 \to A \to B \to P \to 0$$ in $D(H)$-modules, the sequence $$0 \to F(A) \to F(B) \to F(P) \to 0$$ is split as $G$-modules. Therefore, the original sequences are all split, and so $P$ is projective.
\end{proof}

We recall that by \cite[Corollary 5.10]{Balmer1}, $\Spc \stmod(\kk G) \cong \Proj \cohom(G, \kk)$, where $\cohom(G, \kk):=\bigoplus_{i \geq 0} \Ext^i_{\kk G}( \kk, \kk)$ is the cohomology ring of $G$ (recall \exref{fin-group}). 

\bpr{group-dualgroup}
Let $G$, $\kk$, and $\cohom(G,\kk)$ be as above.
\begin{enumerate}
\item The map $f: \Spc \stmod ( \kk G) \to \Spc \stmod( D(\kk G))$ is a homeomorphism, and so $\Spc  \stmod ( D(\kk G)) \cong \Spc \stmod (\kk G) \cong \Proj \cohom(G, \kk)$.
\item Thick ideals of $\stmod ( D(\kk G))$ are in bijection with specialization-closed sets in $\Proj \cohom(G, \kk)$, which are in bijection with thick ideals of $\stmod(\kk G)$, via the maps
$$
\ThickId(\stmod( D(\kk G))) \begin{array}{c} {\Phi_{W}} \atop {\longrightarrow} \\ {\longleftarrow}\atop{\Theta_W} \end{array}  \XX_{\Sp}( \Proj \cohom (G, \kk)) \begin{array}{c} {\Theta_{\kk G}} \atop {\longrightarrow} \\ {\longleftarrow}\atop{\Phi_{\kk G}} \end{array} \ThickId(\stmod(\kk G)).
$$ 
\end{enumerate}
\epr

\begin{proof}
Since $\kk G$ is cocommutative, $\modd( \kk G)$ is braided symmetric. By \leref{group-projectives}, we have $X \cong 0$ in $\StMod(D(H))$ for all $X \in \bK$, and so we are in the situation given of \thref{f-inj-surj-homeo}(3). Additionally, since cohomology rings of groups are finitely generated (for instance by the more general result of \cite{FS1}, in which finite generation of cohomology rings for finite-dimensional cocommutative Hopf algebras in positive characteristic was proven), we know that $\Proj \cohom(G, \kk)$ is a Noetherian topological space. Using Balmer's classification of thick ideals \cite[Theorem 4.10]{Balmer1}, the thick ideals of $\stmod(\kk G)$ are in bijection with specialization-closed sets in $\Spc \stmod(\kk G)$. The rest of the theorem now follows directly as an application of \thref{f-inj-surj-homeo}.
\end{proof}

Now, note that since $\kk[G]^{\cop}$ is a semisimple algebra, $\stmod(\kk[G]^{\cop})$ consists only of the zero object, up to isomorphism, and so $\Spc (\stmod(\kk[G]^{\cop}))$ is the empty set. Thus, the diagram from \reref{hopf-intro} becomes
\begin{center}
\begin{tikzcd}
\Spc \stmod(\kk G) \arrow[rd, "\cong"'] \arrow[r, "\cong"] & {\Proj \cohom(G, \kk)} \arrow[l] & {\Spc \stmod( \kk[G]^{\cop}))= \varnothing} \arrow[ld] \\
                                                           & \Spc \stmod(D(\kk G))               &                                                       
\end{tikzcd}
\end{center}

\subsection{Cosemisimple Hopf algebras}
\label{ss:cosemi}

In fact, we are able to generalize \leref{group-projectives} and \prref{group-dualgroup} from the group algebra case to the case certain finite-dimensional cosemisimple Hopf algebras. Recall that a finite-dimensional Hopf algebra is called {\em{cosemisimple}} if its Hopf dual is semsimple, as an algebra. There has been significant interest in the algebraic properties of cosemisimple Hopf algebras in the past few decades, see e.g. \cite{LR1, LR2, EG1, Chirvasitu, CWW1}. 

We first record the following straightforward lemma. 

\ble{unit-summand-proj}
Let $H$ be a finite-dimensional Hopf algebra such that $\unit_{D(H)}$ is a direct summand of $D(H) \otimes_H \unit_{H}$ as $D(H)$-modules, and let $F: \Mod(D(H)) \to \Mod(H)$ be the forgetful functor. Then $F(P)$ is projective in $\Mod(H)$ if and only if $P$ is projective in $\Mod(D(H))$.
\ele

\begin{proof}
The functor $D(H) \otimes_H -$ is a left adjoint to the forgetful functor $F$. Since $F$ is exact, if $Q$ is a projective $H$-module then 
$$\Hom_H(Q, F(-)) \cong \Hom_{D(H)} ( D(H) \otimes_H Q, -)$$ is an exact functor (recalling that projectives are also injective), and so $D(H) \otimes_H -$ preserves projectivity. Therefore, if $P$ is a $D(H)$-module such that $F(P)$ is projective, then $D(H) \otimes_H F(P)$ is a projective $D(H)$-module. But then, we have $$D(H) \otimes_H F(P) \cong D(H) \otimes_H (\unit_H \otimes_\kk F(P)) \cong (D(H) \otimes_H \unit_H) \otimes_{\kk} P,$$ where the last isomorphism here can be seen from e.g. \cite[Proposition 1.7]{GL1} and the remark following it, which notes that although the proposition is stated for certain universal enveloping algebras, in fact the proof uses only the Hopf algebra structure, and so the result holds for arbitrary Hopf algebras. Note that it holds not just for finite-dimensional modules, but for arbitrary modules, which we need since in this case $P$ may be infinite-dimensional. 

Now, since $\unit_{D(H)}$ is a summand of $D(H) \otimes_H \unit_H$, we have that $P \cong \unit_{D(H)} \otimes_{\kk} P$ is a direct summand of $(D(H) \otimes_H \unit_H) \otimes_{\kk} P$, which is a projective $D(H)$-module, and hence $P$ is projective as well, and the claim is proven.
\end{proof}

Recall that a Hopf algebra (or, more generally, a tensor category) is called {\em{unimodular}} if its spaces of left and right integrals coincide (cf.~ \cite[Section 2.1]{Montgomery1}, \cite[Section 6.5]{EGNO1}). Unimodular Hopf algebras are of particular interest due to their use in constructing Hennings--Kaufman--Radford invariants for 3-manifolds \cite{Hennings1, KR1}. In light of Shimizu's result \cite[Theorem 4.10]{Shimizu1} on unimodular finite tensor categories, if $H$ satisfies the conditions of \leref{unit-summand-proj}-- that is, if $\unit_{D(H)}$ is a direct summand of $D(H) \otimes_H \unit_H$-- then $H$ must be unimodular. The converse is not true; the dual of a finite group algebra is unimodular \cite[Corollary 5.5]{Shimizu1}, but $\unit_{D(\kk[G])}$ is not a direct summand of $D(\kk[G]) \otimes_{\kk G} \unit_{\kk[G]}$ (since $F(\unit_{D(\kk[G])}) = \unit_{\kk[G]}$ is projective and $\unit_{D(\kk[G])}$ is not).

\bco{cosemi-projectives}
Let $H$ be a finite-dimensional unimodular cosemisimple Hopf algebra with Drinfeld double $D(H)$ and forgetful functor $F: \Mod(D(H)) \to \Mod(H)$. Then $F(P)$ is projective as an $H$-module if and only if $P$ is projective as a $D(H)$-module. 
\eco

\begin{proof}
This follows from \leref{unit-summand-proj} and the proof of \cite[Proposition 7.18.15]{EGNO1}. In the course of the proof of the latter, it is shown that if $H$ is unimodular and cosemisimple, then $\unit_{D(H)}$ is a direct summand of $D(H) \otimes_H \unit_{H}$ as $D(H)$-modules (note that here, we are reversing the roles of $H$ and $H^*$ given in their proof). Although this proposition assumes a stronger condition-- that $H$ itself is also semisimple-- this assumption is not used for the part of the proof by which $D(H) \otimes_H \unit_H$ has $\unit_{D(H)}$ as a summand. By \leref{unit-summand-proj}, the corollary follows.
\end{proof}

\bre{kaplansky}
The condition that $H$ is unimodular in \coref{cosemi-projectives} is not too restrictive. It is a long-standing conjecture of Kaplansky \cite{Kaplansky1} that finite-dimensional cosemisimple Hopf algebras are involutory (i.e.~ the square of the antipode is the identity). In view of results of Larson \cite[Corollary 4.2]{Larson}, a weaker form of the Kaplansky conjecture is that all finite-dimensional cosemisimple Hopf algebras are unimodular \cite[Remark 3.9]{AEGN}. This conjecture is still open.
\ere

\coref{cosemi-projectives} and \thref{f-inj-surj-homeo} now immediately imply the following.

\bpr{cosemi-hopf}
Let $H$ be a finite-dimensional unimodular cosemisimple Hopf algebra. Then the map $f: \Spc \stmod(H) \to \Spc \stmod(D(H))$ constructed in Section \ref{ss:mapf} is surjective, and the maps $\Lambda \circ \Psi$ and $\Theta_W \circ \Phi_W$ (as in Section \ref{ss:support}) are each the identity, as maps from the collection of thick ideals of $\stmod(D(H))$ to itself.
\epr

Gelaki has shown \cite[Theorem 1.3.6]{Gelaki1} that every quasitriangular cosemisimple Hopf algebra is unimodular. Hence, again by \coref{cosemi-projectives} and \thref{f-inj-surj-homeo}, we conclude:

\bpr{quasi-cosemi}
Let $H$ be a finite-dimensional quasitriangular cosemisimple Hopf algebra.
\begin{enumerate}
\item The map $f$ constructed in Section \ref{ss:mapf} is a homeomorphism $$\Spc \stmod(H) \xrightarrow{\cong} \Spc \stmod(D(H)),$$ and the maps $\Psi$ and $\Lambda$ as in Section \ref{ss:support} give inverse bijections between the thick ideals of $\stmod(H)$ and $\stmod(D(H))$.
\item If $\Spc \stmod (H)$ is topologically Noetherian, then the $\Phi_W$ and $\Theta_W$ constructed in Section \ref{ss:support} are inverse maps, and so we have the following bijections of thick ideals:
$$
\ThickId(\stmod( D(H))) \begin{array}{c} {\Phi_{W}} \atop {\longrightarrow} \\ {\longleftarrow}\atop{\Theta_W} \end{array}  \XX_{\Sp}( \Spc( \stmod(H))) \begin{array}{c} {\Theta_{H}} \atop {\longrightarrow} \\ {\longleftarrow}\atop{\Phi_{H}} \end{array} \ThickId(\stmod(H)).
$$ 
\end{enumerate}
\epr

Of course, if $H$ itself is also semisimple, then \prref{cosemi-hopf} and \prref{quasi-cosemi} are not particularly illuminating, since this implies that $D(H)$ is also semisimple, and then the Balmer spectra of $\stmod(H)$ and $\stmod(D(H))$ are both $\varnothing$. It is a classical theorem of Larson--Radford \cite{LR1} that in characteristic 0, all cosemisimple finite-dimensional Hopf algebras are also semisimple. Hence, \prref{cosemi-hopf} and \prref{quasi-cosemi} only provide interesting examples in positive characteristic.

\subsection{Benson--Witherspoon smash coproduct Hopf algebras}
\label{ss:bwhopf}

We will now consider the Benson--Witherspoon smash coproducts which were originally studied in \cite{BW1}, with generalizations studied in \cite{MVW1} and \cite{PW1}; their Balmer spectra and thick ideals were classified in \cite{NVY1}. We recall the general construction of these algebras. Let $G$ and $L$ be finite groups, such that $L$ acts on $G$ by group automorphisms, and let $\kk$ be an algebraically closed field of characteristic dividing the order of $G$. We then define $H_{G,L}$ to be the Hopf algebra dual of the smash product $\kk [G] \# \kk L$, where $\kk [G]$ is the coordinate ring of $G$, and $\kk L$ is the group algebra of $L$. 

As an algebra, $H_{G,L}$ is isomorphic to $\kk G \otimes \kk [L]$. We will denote by $\{ p_x : x \in L\}$ the standard dual basis for $\kk [L]$, as in Section \ref{ss:groupalg}. Denote by $e$ the identity element of $L$. The additional Hopf algebra structures of comultiplication, counit, and antipode on $A$ are defined by
\begin{align*}
\Delta(g \otimes p_x) &= \sum_{y \in L} ( g \otimes p_y ) \otimes (y^{-1} . g \otimes p_{y^{-1} x } ),\\
\epsilon( g \otimes p_x ) &= \delta_{x,1},\\
S(g \otimes p_x) &= x^{-1} . (g^{-1} ) \otimes p_{x^{-1}},
\end{align*}
for all $g \in G$ and $x \in L$.

Since as an algebra $H_{G,L} \cong \kk G \otimes \kk [L]$, an $H_{G,L}$-module is the same as a $G$-module with an $L$-grading, such that the action of $G$ preserves the $L$-grading. That is, every $H_{G,L}$-module $M$ may be decomposed 
$$M \cong \bigoplus_{x \in L} M_x \otimes \kk_x$$
where $M_x$ is a $G$-module, and $\kk_x$ is the 1-dimensional $\kk[L]$-module on which $p_x$ acts as the identity, and $p_y$ acts as 0 for $y \not = x$ (in other words, the $\kk[L]$-module corresponding to a $L$-graded vector space of one dimension where every element is homogeneous of degree $x$). The $H_{G,L}$-action on the component $M_x \otimes \kk_x$ is defined by letting $\kk G$ act on the first tensorand, and $\kk[L]$ act on the second. 

Using the definition of the coproduct on $H_{G,L}$, Benson and Witherspoon \cite[Theorem 2.1]{BW1} compute the formula for the tensor product of $H_{G,L}$-modules:
$$(M_x \otimes \kk_x) \otimes (N_y \otimes \kk_y) = ( M_x \otimes {^x N_y} ) \otimes \kk_{xy},$$
for any $\kk G$-modules $M_x$ and $N_y$, and for all $x, y \in L$, where the module $^x N_y$ is defined as the twist of the module $N_y$ by the action of $x$. Namely, this is the $\kk G$-module which is equal to $N_y$ as a vector space, and if we write $g \cdot v$ for the action of $G$ on the original module $N_y$, then the new action $*$ of $G$ on $^x N_y$ is defined $g * v = (x^{-1}g) \cdot v.$

\bpr{bw-drinf}
Let $H_{G,L}$ the Benson--Witherspoon smash coproduct Hopf algebra as defined above, $\bC$ the category $\modd(H_{G,L})$, and $\drin(\bC)$ the category $\modd( D(H_{G,L}))$ for the Drinfeld double $D(H_{G,L})$ of $H_{G,L}$. 
\begin{enumerate}
\item The continuous map $f: \Spc \st(\bC) \to \Spc \st(\drin(\bC))$ constructed in Section \ref{ss:mapf} is injective.
\item The map $\Psi \circ \Lambda$ constructed in Section \ref{ss:support} is equal to the identity, as a map $\ThickId(\st(\bC)) \to \ThickId(\st(\bC))$.
\item The map $\Phi_W \circ \Theta_W$ constructed in Section \ref{ss:support} is equal to the identity, as a map $\mc{X}_{\Sp} ( \Spc \st(\bC)) \to \mc{X}_{\Sp} (\Spc \st(\bC))$. 
\end{enumerate}
\epr

\bre{not-follows-thm}
We note that if $\bC$ was braided, then \prref{bw-drinf} would follow directly from \thref{f-inj-surj-homeo}. However, in general, $H_{G,L}$ is not a quasitriangular Hopf algebra, i.e.~ the category of $H_{G,L}$-modules is not braided. 
\ere

\prref{bw-drinf} will be proven by first showing the following intermediary lemma.

\ble{det-by-cent}
Suppose $\bI$ and $\bJ$ are thick ideals of $\st(\bC)$ such that
$$\{ X \in \st(\drin(\bC)) : \overline{F} (X) \in \bI \} = \{ X  \in \st(\drin(\bC)): \overline{F} (X) \in \bJ \}.$$
Then $\bI = \bJ$. In particular, if $M$ is an object of $\st(\drin(\bC))$, then there exists an object $\hat M$ which is in the image of $\overline{F}$, and given any thick ideal $\bI$, the object $M$ is in $\bI$ if and only if $\hat M$ is in $\bI$. 
\ele

\begin{proof}
Suppose $\bI$ and $\bJ$ are thick ideals satisfying the condition above. Since $\bI$ and $\bJ$ are thick, it is enough to show that the indecomposable objects in $\bI$ are equal to the indecomposable objects in $\bJ$. Suppose $M_x \otimes \kk_x$ is an object in $\bI$. Then the module
$$(M_x \otimes \kk_x) \otimes ( \kk \otimes \kk_{x^{-1}}) \cong M_x \otimes \kk_{e}$$ 
is in $\bI$. We also then have
$$( \kk \otimes \kk_y) \otimes (M_x \otimes \kk_{e} ) \otimes ( \kk \otimes \kk_{y^{-1}})  \cong \; {^y M_x} \otimes \kk_{e}$$
is an object of $\bI$ as well. The ideal $\bI$ then contains the direct sum
$$ \hat M := \bigoplus_{y \in H} \;^y M_x \otimes \kk_{e}.$$
We claim that $\hat M$ is in the image of $\overline{F}$; in other words, $\hat M$ has a half-braiding which allows it to be lifted to the Drinfeld center. To see this, consider an $H_{G,L}$-module $N_z \otimes \kk_z$. We observe that
\begin{align*}
\hat M \otimes (N_z \otimes \kk_z) &\cong \bigoplus_{y \in L} ({^y M_x} \otimes N_z) \otimes \kk_z, \\
(N_z \otimes \kk_z) \otimes \hat M &\cong \bigoplus_{y \in L} ( N_z \otimes {^{zy}M_x} ) \otimes \kk_z. \\
\end{align*}
Since $\kk G$ is itself cocommutative (and thus $^yM_x \otimes N_z \cong N_z \otimes {^yM_x}$ in a natural way), this formula can be used to observe a natural isomorphism $\hat M \otimes - \cong - \otimes \hat M$. This isomorphism satisfies the half-braiding condition, and so $\hat M$ is in the image of $\overline{F}$.

Since $\bI$ and $\bJ$ are assumed to agree on their intersections with the image of $\overline{F}$, we can conclude that $\hat M$ is in $\bJ$ as well. But then its summand $M_x \otimes \kk_{e}$, and hence $$(M_x \otimes \kk_{e} ) \otimes ( \kk \otimes \kk_{x}) \cong M_x \otimes \kk_x,$$ is also an object of $\bJ$. Note that we have proven generally that $M_x \otimes \kk_x$ is in any thick ideal if and only if $\hat M$, as constructed above, is in that ideal. Thus, the objects of $\bI$ are a subset of the objects of $\bJ$, and by symmetry the ideals are equal.
\end{proof}

We can now prove \prref{bw-drinf}, as a consequence of \leref{det-by-cent}:

\begin{proof}
The map $f$ is defined by $$f(\bP) = \{ X \in \st(\drin(\bC)): \overline{F} (X) \in \bP \}$$ for a given prime ideal $\bP$ in $\Spc \st(\bC)$. But \leref{det-by-cent} has shown that if $\bP$ and $\bQ$ are two prime ideals with $f(\bP) = f(\bQ)$, then since $\bP$ and $\bQ$ are more generally examples of thick ideals, we have $\bP = \bQ$. Hence, $f$ is injective, showing (1).

For (2), let $S$ be an arbitrary specialization-closed set in $\Spc \st(\bC)$, in other words, a (possibly infinite) union $S= \bigcup_{i\in I} S_i$ where each $S_i$ is a closed set. Recall that by construction, it is automatic that $\Phi_W (\Theta_W(S)) \subseteq S$ (the details are included above in the proof of \thref{f-inj-surj-homeo}). 

To show the opposite containment, we note that by the classification of thick ideals and Balmer spectrum of $\st(\bC)$ as given in \cite{NVY1}, $\Spc \st(\bC)$ is a Noetherian topological space. We claim that this implies that every closed set in $\Spc \st(\bC)$ has the form $V_{\bC}(M)$, for some object $M$ of $\st(\bC)$, just as in the commutative setting \thref{balmer-top}, using \leref{det-by-cent} as a substitute for the commutativity of the tensor product. Let $V_{\bC}(\mc{T})$ be an arbitrary closed set in $\Spc \st(\bC)$, for some collection $\mc{T}$ of objects in $\st(\bC)$. Then the complement of $V_{\bC}(\mc{T})$ is by definition
$$U_{\bC}(\mc{T}) = \{\bP \in \Spc \st(\bC) : \bP \cap \mc{T} \not = \varnothing\},$$
and has an open cover
$$U_{\bC}(\mc{T}) = \bigcup_{A \in \mc{T}} U_{\bC}(A) = \bigcup_{A \in \mc{T}} \{ \bP \in \Spc \st(\bC): A \in \bP\}.$$
By Noetherianity, this set is compact, and hence has a finite subcover
$$U_{\bC}(\mc{T}) = \bigcup_{A \in \mc{T}'} U_{\bC}(A), $$
where $\mc{T}' \subseteq \mc{T}$ is some finite collection of objects. Enumerate the objects of $\mc{T}'$ by $A_1,..., A_n$. Choose $\hat A_1,..., \hat A_n$ as constructed in \leref{det-by-cent}: they are in the image of $\overline{F}$, and for any thick ideal $\bI$, we have $A_j \in \bI$ if and only if $\hat A_j \in \bI$. Using this property, it is clear that $V_{\bC}(A_j)=V_{\bC}(\hat A_j)$ for all $j$. Now we claim that
$$U_{\bC}(\mc{T}) = U_{\bC}(A_1) \cup... \cup U_{\bC}(A_n) = U_{\bC}(\hat A_1) \cup... \cup U_{\bC}(\hat A_n) = U_{\bC}(\hat A_1 \otimes... \otimes \hat A_n).$$ The last equality (more specifically, the containment $\supseteq$) uses the fact that each $\hat A_j$ is in the image of $\overline{F}$, and hence commutes with all objects of $\st(\bC)$ up to isomorphism, since this implies that
$$\hat A_1 \otimes... \otimes \hat A_n \in \bP \Rightarrow A_j \in \bP \text{ for some }j.$$
Our claim is now shown: every closed set in $\Spc \st(\bC)$ is of the form $V_{\bC}(A)$ for some object $A$. 

In particular, each of the closed sets $S_i$, for $i \in I$, can be written as $V_{\bC} ( M_i)$ for some object $M_i \in \st(\bC)$. As above, we can replace $M_i$ by $\hat M_i$, which is is in the image of $\overline{F}$, i.e. we can pick an object $X_i$ in $\st(\drin(\bC))$ with $\overline{F}(X_i) = \hat M_i$. Since $$W(X_i)= V_{\bC} ( \overline{F}(X_i)) = V_{\bC} ( \hat M_i) = V_{\bC} (M_i) = S_i \subseteq S,$$ we have $X_i \in \Theta_W (S)$ by definition. Hence, we now have
\begin{align*}
\Phi_W ( \Theta_W (S)) &\supseteq \bigcup_{i \in I} W(X_i)\\
&= \bigcup_{i \in I} S_i \\
&= S.
\end{align*}
Since we have both containments, we can conclude that $\Phi_W (\Theta_W (S)) = S$ for any specialization-closed set $S$ in $\Spc \st(\bC)$.
\end{proof} 

We also note that if $p$ does not divide the order of $L$, then we can apply the results of the previous section to obtain:

\bth{cosemi-bw}
Let $G$, $L$, $\kk$, $H_{G,L}$, $\bC= \modd(H_{G,L})$, and $\drin(\bC)= \modd(D(H_{G,L}))$ be as above, and assume additionally that $p$ does not divide the order of $L$. Then we have the following.
\begin{enumerate}
\item The map $f$ constructed in Section \ref{ss:mapf} is a homeomorphism $$\Spc \st(\bC) \xrightarrow{\cong} \Spc \st(\drin(\bC)).$$
\item The maps $\Phi_W$ and $\Theta_W$ constructed in Section \ref{ss:support} are mutually inverse, and so we have the following bijections of thick ideals:
$$
\ThickId(\stmod( D(H_{G,L}))) \begin{array}{c} {\Phi_{W}} \atop {\longrightarrow} \\ {\longleftarrow}\atop{\Theta_W} \end{array}  \XX_{\Sp}( \Spc( \stmod(H_{G,L}))) \begin{array}{c} {\Theta_{H_{G,L}}} \atop {\longrightarrow} \\ {\longleftarrow}\atop{\Phi_{H_{G,L}}} \end{array} \ThickId(\stmod(H_{G,L})).
$$ 
\end{enumerate}
\eth

\begin{proof}
First, note that $H_{G,L}$ is cosemisimple: its dual is the smash product $\kk[G] \# \kk L$. Since $p$ does not divide the order of $L$, the group algebra $\kk L$ is semisimple, and by \cite[Theorem 6]{CF1}, as the smash product of two semisimple algebras, $\kk[G] \# \kk L$ is semisimple as well.

Next, we claim that $H_{G,L}$ is unimodular. This can be observed directly, by noting that the element 
$$h:=\left ( \sum_{g \in G} g \right ) \otimes p_1$$ 
is both a left and a right integral in $H_{G,L}$.

By application of \prref{bw-drinf} and \prref{cosemi-hopf}, $f$ is bijective and the maps $\Phi_W$ and $\Theta_W$ are inverse bijections. To conclude the proof, we must just prove that $f$ is closed, and hence a homeomorphism. This follows similarly to the proof of \thref{f-inj-surj-homeo}(3a), except that we must again use \leref{det-by-cent} as a substitute for commutativity of the tensor product. Let $V_{\bC}(M)$ an arbitrary closed set, and, just as before, we may assume (by replacing $M$ with $\hat M$ as in \leref{det-by-cent} if need be) that $M$ is in the image of $\overline{F}$, and so we can pick $X \in \drin(\st(\bC))$ with $\overline{F}(X)= M$. We now have
\begin{align*}
f(V_{\bC}(M))&=\{ f(\bP): \bP \in \Spc \st(\bC), M \not \in \bP\}\\
&=\{ \bQ \in \Spc \st(\drin(\bC)): X \not \in \bQ \}\\
&= V_{\drin}(X).
\end{align*}
The second equality follows from the fact that $f$ is bijective. Hence, $f$ is closed, and the theorem is complete. 
\end{proof}

\end{document}